\documentclass[12pt]{amsart}

\usepackage{longtable}
\usepackage{float}
\restylefloat{table}

\newcommand{\hide}[1]{}

\usepackage{amssymb}
\usepackage{amsbsy}
\usepackage{amscd}
\usepackage{amsmath}
\usepackage{amsthm}
\usepackage{bbold}
\usepackage{mathrsfs}
\usepackage{url}
\usepackage{graphicx}

\numberwithin{equation}{section}

\theoremstyle{plain}
\newtheorem{thm}{Theorem}[section]
\newtheorem{prop}[thm]{Proposition}

\newtheorem{thm-defi}[thm]{Theorem/Definition}
\newtheorem{cor}[thm]{Corollary}
\newtheorem{lem}[thm]{Lemma}

\theoremstyle{definition}

\newtheorem{question}[thm]{Question}
\newtheorem{example}[thm]{Example}

\newcommand{\A}{{\mathcal A}}

\newcommand{\CC}{{\mathbb C}}

\newcommand{\E}{{\mathcal E}}

\newcommand{\G}{{\mathcal G}}

\newcommand{\QQ}{{\mathbb Q}}

\newcommand{\RR}{{\mathbb R}}

\newcommand{\T}{{\mathcal T}}

\newcommand{\U}{{\mathcal U}}

\newcommand{\Y}{{\mathcal Y}}

\newcommand{\ZZ}{{\mathbb Z}}
\renewcommand{\P}{{\mathcal P}}
\newcommand{\PP}{{\mathbb P}}

\newcommand{\Spin}{{\rm Spin}}

\newcommand{\IsomRightArrow}{\stackrel{\cong}{\rightarrow}}

\newcommand{\RightArrowOf}[1]{\stackrel{#1}{\rightarrow}}

\newcommand{\LongRightArrowOf}[1]{\stackrel{#1}{\longrightarrow}}

\newcommand{\StructureSheaf}[1]{{\mathcal O}_{#1}}
\newcommand{\EndProof}{\hfill  $\Box$}

\newcommand{\Gal}{{\rm Gal}}

\newcommand{\rank}{{\rm rank}}

\newcommand{\Pic}{{\rm Pic}}

\newcommand{\Sym}{{\rm Sym}}
\newcommand{\Ext}{{\rm Ext}}

\newcommand{\Hom}{{\rm Hom}}
\newcommand{\Aut}{{\rm Aut}}
\newcommand{\End}{{\rm End}}

\newcommand{\RHom}{{\mathcal RH}om}
\newcommand{\SheafHom}{{\mathcal H}om}

\newcommand{\SheafExt}{{\mathcal E}xt}

\newcommand{\Ideal}[1]{{\mathcal I}_{#1}}

\newcommand{\Contract}{\rfloor}

\newcommand{\Choose}[2]
{\left(\!\!\begin{array}{c}#1\\#2\end{array}\!\!\right)}

\renewcommand{\span}{{\rm span}}

\input xy
\xyoption{all}

\begin{document}
\title[Secant sheaves  
and Weil classes on abelian varieties]
{Secant sheaves and Weil classes on abelian varieties}
\author{Eyal Markman}
\address{Department of Mathematics and Statistics, 
University of Massachusetts, Amherst, MA 01003, USA. Email: markman@umass.edu}

\date{\today}

\begin{abstract}
Let $K$ be a CM-field, i.e., a totally complex quadratic extension of a totally real field $F$. Let $X$ be an abelian variety admitting an algebra embedding $F\rightarrow \End_\QQ(X)$, and let 
$\hat{X}$ be the dual abelian variety. We construct an embedding $\eta:K\rightarrow \End_\QQ(X\times\hat{X})$ associated to a choice of a polarization $\Theta$ in $\wedge^2_FH^1(X,\QQ)$ and an element $q\in F$, such that $K=F(\sqrt{-q})$. 
We get the $[K:\QQ]$-dimensional subspace $HW(X\times\hat{X},\eta)$ 
of Hodge Weil classes in $H^{\frac{d}{2},\frac{d}{2}}(X\times\hat{X},\QQ)$, where $d:=4\dim(X)/[K:\QQ]$. 

The space $V_\CC\!:=\!H^1(X\times\hat{X},\CC)$ admits a natural symmetric bilinear pairing and 
the even cohomology $S^+_\CC\!:=\!H^{ev}(X,\CC)$ is the half-spin representation of $\Spin(V_\CC)$. Hence, $\PP(S^+_\CC)$ contains 
a component of the Grassmannian of maximal isotropic subspaces of $V_\CC$ known as 
the even spinorial variety. 
We associate to $(\Theta,q)$ a $2^{[F:\QQ]}$-dimensional subspace $B$ of $S^+_\QQ$ such that $\PP(B)$ is secant to the spinorial variety. Associated to two coherent sheaves $F_1$ and $F_2$ on $X$ with Chern characters in $B$ we obtain the object $E:=\Phi(F_1\boxtimes F_2^\vee)$ in the derived category $D^b(X\!\times\!\hat{X})$, where $\Phi\!:\!D^b(X\!\times\! X)\!\rightarrow \! D^b(X\!\times\!\hat{X})$ is Orlov's equivalence.
The flat deformations of the normalized Chern character $\kappa(E):=ch(E)\exp\left(-\frac{c_1(E)}{\rank(E)}\right)$ of $E$ remain of Hodge type under every deformation of $(X\times\hat{X},\eta)$ as an abelian variety  of Weil type $(A',\eta')$. The algebraicity of the Weil classes of every deformation $(A',\eta')$ would thus follow if E is semiregular in the appropriate sense.

When $F=\QQ$, so that $K$ is an imaginary quadratic number field,
the above 
construction was combined with the Semi-regularity theorem  
to prove the algebraicity of the Weil classes on abelian sixfolds of split Weil type. The algebraicity of the Weil classes on all abelian fourfold of Weil type follows. The Hodge conjecture for abelian varieties of dimension $\leq 5$ is known to follow from the latter result.
\end{abstract}

\maketitle

\setcounter{tocdepth}{1}
\tableofcontents

%
\section{Introduction}
\label{sec-introduction}

%
\subsection{Abelian varieties of Weil type}
A {\em CM-field} $K$ is a quadratic extension of a number field $F$, such that all embeddings of $F$ in $\CC$ are real and none of the embeddings of $K$ in $\CC$ are real. Set $e:=[K:\QQ]$. Let $\Sigma$ be the set of all embeddings $\sigma:K\rightarrow \CC$. The cardinality of $\Sigma$ is $e$. 
Let $\iota$ be the involution in $\Gal(K/F)$. Then $\sigma\circ\iota=\bar{\sigma}$, for all $\sigma\in\Sigma$, where $\bar{\sigma}$ is the complex conjugate embedding.
Let $A$ be a complex (projective) abelian variety and $\eta:K\rightarrow \End_\QQ(A):=\End(A)\otimes_\ZZ\QQ$ an algebra embedding. Note that $\eta$ is equivalent to an embedding
$\eta:K\rightarrow \End_{Hdg}(H^1(A,\QQ))$. 
Let $H^1_\sigma(A,\CC)$ be the subspace of $H^1(A,\CC)$ on which $\eta(K)$ acts via the character $\sigma\in\Sigma$ and set $H^{1,0}_\sigma(A):=H^{1,0}(A)\cap H^1_\sigma(A,\CC)$.
Define $H^{0,1}_\sigma(A)$ analogously. Set $d:=\dim_K H^1(A,\QQ)$.
The pair 
$(A,\eta)$ is said to be of {\em Weil type}, if $\dim H^{1,0}_\sigma(A)=\dim H^{0,1}_{\sigma}(A)=\frac{d}{2}$, for all $\sigma\in\Sigma$.
In that case the subspace\footnote{
We have the decomposition $K\otimes_\QQ\CC\cong \oplus_{\sigma\in\Sigma}\CC$ yielding  
$
H^1(A,\CC)=H^1(A,\QQ)\otimes_\QQ\CC =H^1(A,\QQ)\otimes_K (K\otimes_\QQ\CC)\cong \oplus_{\sigma\in\Sigma} H^1_\sigma(A,\CC).
$
The subspace $\oplus_{\sigma\in\Sigma} \wedge^d\!\!H^1_\sigma(A,\CC)$ of $H^d(A,\CC)$ is defined over $\QQ$ and corresponds to the subspace 
$\wedge^d_KH^1(A,\QQ)$ of $\wedge^d\!H^1(A,\QQ)$.
} 
$HW(A,\eta):=\wedge^d_KH^1(A,\QQ)$ of $\wedge^dH^1(A,\QQ)$ is an $e$-dimensional subspace 
of $H^{\frac{d}{2},\frac{d}{2}}(A,\QQ)$, which is a $1$-dimensional $K$-vector space \cite[Prop. 4.4]{deligne-milne}. 
Classes in $HW(A,\eta)$ are called {\em Weil classes}. Note that both $e$ and $d$ are necessarily even, and so $\dim_\CC(A)=2n$, where 
$
n:=de/4
$
is an integer.
\[
e:=[K:\QQ], \ \ \ d:=\dim_KH^1(A,\QQ), \ \ \ n:=\dim_\CC(A)/2=de/4.
\]

A {\em polarized abelian variety of Weil type} is a triple $(A,\eta,h)$, where\footnote{
Let $\hat{\Sigma}$ be the set of embeddings of $F$ into $\RR$. We get the decomposition $H^1(A,\RR)=\oplus_{\hat{\sigma}\in\hat{\Sigma}}H^1_{\hat{\sigma}}(A,\RR)$. The subspace $\oplus_{\hat{\sigma}\in\hat{\Sigma}}\wedge^2H^1_{\hat{\sigma}}(A,\RR)$ is defined over $\QQ$ and corresponds to the subspace $\wedge^2_FH^1(A,\QQ)$.
} 
$h\in \wedge^2_F H^1(A,\QQ)\cap H^{1,1}(A,\QQ)$ is the class of a polarization satisfying $h(\eta_t(x),y)=h(x,\eta_{\iota(t)}(y))$, for all $t\in K$, $x,y\in H_1(A,\QQ)$.
Such a polarization yields a non-degenerate $K$-valued hermitian form $H:H_1(A,\QQ)\times H_1(A,\QQ)\rightarrow K$ on the $K$ vector space
$H_1(X,\QQ)$ and $(A,\eta,h)$ is said to be of {\em split Weil type}, if $H$ has an isotropic subspace of half the dimension.
The moduli space of polarized abelian varieties of Weil type is $ed^2/8$-dimensional. Assume that $d\geq 4$. For a generic triple $(A,\eta,h)$ in moduli, 
the rank of the Neron-Severi group $H^{1,1}(A,\ZZ)$ is $e/2$, 
the Hodge ring of $A$ is generated by $H^{1,1}(A,\QQ)$ and $HW(A,\eta)$, and $HW(A,\eta)$ intersects trivially the subalgebra generated by  
$H^{1,1}(A,\QQ)$. 
In particular, for a generic triple $(A,\eta,h)$ in moduli, 
\begin{equation}
\label{eq-Hodge-classes-of-degree-d-on-the-generic-A}
H^{\frac{d}{2},\frac{d}{2}}(A,\QQ)=Im[\Sym^{d/2}\left(H^{1,1}(A,\QQ)\right)]
\oplus HW(A,\eta),
\end{equation}
where the first summand is the image of $\Sym^{d/2}\left(H^{1,1}(A,\QQ)\right)$ in $H^{\frac{d}{2},\frac{d}{2}}(A,\QQ)$.
The Hodge conjecture suggests the following.

\begin{question}
Does $HW(A,\eta)$ consist of algebraic classes?
\end{question}

The question of algebraicity of the Weil classes was considered by some as a first test case for the Hodge conjecture. 
The following historical note is part of \cite[footnote 14]{milne-Tate}: 
\begin{small}
{\em Mumford and Tate tried to prove the Hodge conjecture for abelian varieties by showing that the $\QQ$-algebra of rational $(p,p)$ classes is generated by those of type $(1,1)$, for which the conjecture was known, but Mumford found a counterexample to this. When Tate told Weil of the example, he remarked that it is a special case of a slightly more generic example, namely a $4$-dimensional family of examples, and then said
``As you and Mumford seem to believe Hodge's conjecture, it is up to you to exhibit algebraic cycles corresponding to these abnormal classes. I shall rather attempt to show there is no such cycle'' (Letter from Tate to Serre, February 2, 1965.) 
\dots
}
\end{small}
See  \cite{zharkov} for a computer aided attempt to find counter examples to the algebraicity of the Weil classes via degeneration and tropical algebraic geometry.

Weil's 1977 paper \cite{weil} contains the construction mentioned in the quote above. He constructs the moduli spaces of abelian varieties of Weil type and their Weil classes for $K$ an imaginary quadratic number field. The more general construction for CM-fields can be found in \cite{deligne-milne}.

\begin{thm} \cite[Theorem 1.5.1]{markman-sixfolds}
\label{main-thm}
The Weil classes for abelian fourfolds of Weil type and abelian sixfolds of split Weil type with complex multiplication by a quadratic imaginary number field $K$
are algebraic.
\end{thm} 

A sketch of the proof of Theorem \ref{main-thm} is given in Section \ref{sec-semiregular-sheaf-over-sixfold}.
The above result was proved earlier 
by Schoen in case $K=\QQ(\sqrt{-3})$ \cite{schoen} and by Koike in case $K=\QQ(\sqrt{-1})$ \cite{koike}. The algebraicity for abelian fourfolds of split Weil type was proven earlier for $K=\QQ(\sqrt{-3})$ in
\cite{schoen1}, for $K=\QQ(\sqrt{-1})$ in \cite{van-Geemen}, and for all quadratic imaginary number fields in \cite{markman-generalized-kummers}. An alternative proof for abelian fourfolds of split Weil type was obtained recently by Floccari and Fu \cite{FF}. The proofs in \cite{FF,markman-generalized-kummers} use the geometry of hyper-K\"{a}hler varieties.

As we shall see, the proof of Theorem \ref{main-thm} relies on developments which were not available at the time Weil made his remark. Among these developments are the works of Mukai, Polishchuk, and Orlov on equivalences of derived categories of abelian varieties and the works of Bloch, Buchweitz-Flenner, and Pridham on the Semi-regularity theorem \cite{bloch,buchweitz-flenner,mukai-duality,orlov-abelian-varieties,polishchuk,pridham}. Furthermore, the strategy of the proof is a generalization of the one in \cite{markman-generalized-kummers}, which in turn is inspired by O'Grady's observation that the third intermediate Jacobians of projective hyper-K\"{a}hler varieties of Kummer type form complete $4$-dimensional families of abelian fourfolds of split Weil type \cite{ogrady}.
 
 The Hodge conjecture is known for projective varieties of dimension $\leq 3$.
The Hodge ring of abelian fourfolds is generated by divisor classes and Weil classes for complex multiplication by possibly more than one imaginary quadratic number field, by
work of Moonen and Zarhin \cite{Moonen-Zarhin-Weil-Hodge-Tate-classes,moonen-zarhin-low-dimension} 
 combined with a result of Ram\'{o}n Mari in the case of products of abelian surfaces
\cite{ramon-mari}. The Hodge ring for simple abelian varieties of prime dimension is generated by divisor classes, by a result of Tankeev \cite{tankeev}.
If $X$ is a non-simple abelian variety of dimension $5$, then the Hodge ring of $X$ is generated by divisor classes and pull backs of Weil classes from quotient abelian fourfolds, by \cite[Theorem 0.2]{moonen-zarhin-low-dimension}.
Combining these results with Theorem \ref{main-thm} we get:

\begin{cor} 
The Hodge conjecture holds for abelian varieties of dimension $\leq 5$.
\end{cor}

We will present in this paper a general strategy for proving the algebraicity of the Weil classes on abelian varieties of split Weil type but implement it fully only for dimension $\leq 6$ and $K$ imaginary quadratic. Following is a class of abelian varieties for which the Hodge conjecture would follow, if the algebraicity of the Weil classes of split type would be proved.
An abelian variety $A$ is of {\em CM-type}, if it admits an algebra embedding $\eta:K\rightarrow \End_\QQ(A)$, such that $H^1(A,\QQ)$ is a $1$-dimensional $K$-vector space. 
In this case $[K:\QQ]=2\dim(A)$ and we have the decomposition $H^1(A,\CC)=\oplus_{\sigma\in \Sigma}H^1_\sigma(A,\CC)$ by $1$-dimensional subspaces. If  
$H^{1,0}_\sigma(A)\neq(0)$, then $H^{0,1}_{\bar{\sigma}}(A)\neq(0)$. Hence, $H^{1,0}(A)=\oplus_{\sigma\in T}H^{1,0}_\sigma(A)$, where $T\subset \Sigma$ consists of precisely one embedding of $K$ out of each pair of complex conjugate embeddings. Such a pair $(K,T)$ is called a {\em CM-type}.
Every simple abelian variety of CM-type is of the form 
$\CC^g/T(\mathfrak{a})$, where $2g=[K:\QQ]$, $\mathfrak{a}$ is an ideal in the subring of $K$ of algebraic integers, and $T(\mathfrak{a})$ is the rank $2g$ lattice in $\CC^g$ consisting of $\{(\sigma_1(k), \dots, \sigma_g(k)) \ : \ k\in \mathfrak{a}\}$, where $T=\{\sigma_i \ : \ 1\leq i\leq g\}$. Abelian varieties of CM-type are rigid. They correspond to special points in Shimura varieties and play a central role in the theory of Shimura varieties \cite{milne-CM,kerr}.
Andr\'{e} reduced the Hodge conjecture for abelian varieties of CM-type to the question of algebraicity of the Weil classes on abelian varieties of split Weil type.

\begin{thm} \cite{andre}
\label{thm-andre}
Let $A$ be a complex abelian variety of CM-type. There exist abelian varieties $A_i$ of split Weil type and homomorphisms $f_i:A\rightarrow A_i$, such that every Hodge class $t$ on $A$ can be written as a sum $t=\sum f_i^*(t_i)$ with $t_i$ a Weil class on $A_i$.
\end{thm}

The Hodge conjecture for CM abelian varieties implies Grothendieck's standard conjecture for all abelian varieties, the Tate conjecture for all abelian varieties over finite fields, and it makes it possible to implement Milne's ``program to extend Deligne's theory of absolute Hodge classes to characteristic $p$, thereby obtaining a good theory of abelian motives in mixed characteristic''
\cite[Footnote 3]{milne-abelian-motives}.
%
\subsection{Organization of the paper}
A general strategy for proving the algebraicity of the Weil classes in $HW(A,\eta)$, $\eta:K\rightarrow \End_\QQ(A)$, for imaginary quadratic number fields $K$, was developed in \cite{markman-sixfolds} and implemented for abelian fourfolds of Weil type and abelian sixfolds of split Weil type. In sections \ref{sec-semiregularity} to \ref{sec-BB}
we present\footnote{Details appear in \cite{markman-CM} due to page limitations for this contribution.} the natural generalization of this strategy for $K$ a CM-field. In section \ref{sec-semiregular-sheaf-over-sixfold} we survey the implementation of the strategy in the case of imaginary quadratic $K$ carried out in \cite{markman-sixfolds}.

In Section \ref{sec-semiregularity} we recall the Semi-regularity theorem of Buchweitz and Flenner.
In Section \ref{sec-spin-groups} we recall the results of Mukai, Polishchuk, and Orlov relating (1) the action on $H^*(X,\ZZ)$ of the group of autoequivalences of the derived category $D^b(X)$ of an abelian variety $X$, and (2) the spin representation of the group $\Spin(V)$ of the lattice $V=H^1(X,\ZZ)\oplus H^1(X,\ZZ)^*$ endowed with its natural symmetric bilinear pairing. 
In Section \ref{sec-strategy} we formulate a general strategy for proving the algebraicity of Weil classes guided by the Semi-regularity theorem and Chevalley's theory of pure spinors. 

Set $V_\bullet:=V\otimes_\ZZ\bullet$, $\bullet=\ZZ$ or a field.
In Section \ref{sec-pure-spinor} we review how the theory of pure spinors relates classes in the spin representation $H^*(X,\QQ)$ of $\Spin(V)$ to maximal isotropic subspaces of $V_K$. When $X$ admits an algebra embedding  $\hat{\eta}:F\rightarrow \End_\QQ(X)$ of the totally real subfield $F$ of $K$,
we associate to a maximal isotropic subspace $W$ of $V_K$ a collection of maximal isotropic subspaces $W_T$ in $V_\CC$, one for each CM-type $T$ of $K$.
We introduce the subspace $B$ of $H^*(X,\QQ)$ spanned by the subset of the spinorial variety consisting of the pure spinors of $\{W_T\}$.  By definition, $\PP(B)$ is secant to the spinorial variety. 
In Section \ref{sec-complex-multiplication} we observe that if $B$ is spanned by Hodge classes, then a natural embedding $\eta:K\rightarrow V_\QQ$ associated to $W$ endows $X\times\hat{X}$ with the structure of an abelian variety of Weil type. 

Let $\Spin(V_\QQ)_B$ be the subgroup of $\Spin(V_\QQ)$ fixing every class in $B$.
In Section \ref{sec-A} we compute the subalgebra $\A$ of $\Spin(V_\QQ)_B$-invariant classes in $H^*(X\times\hat{X},\QQ)$. We show that $\A$ is generated by its graded summand $\A^2\subset H^{1,1}(X\times\hat{X},\QQ)$ and by the subspace of Weil classes $HW(X\times\hat{X},\eta).$ 
Let $K_-\subset K$ be the $-1$-eigenspace of the involution $\iota\in\Gal(K/F)$. We construct an isomorphism $\Xi:K_-\rightarrow \A^2$ and associate to $t\in K_-$ a 
$K$-valued hermitian form $H_t$ on the $K$-vector space $(V_\QQ,\eta)$. Let $m:\Spin(V_\RR)\rightarrow GL(H^*(X,\RR))$ and $\rho:\Spin(V_\RR)\rightarrow SO_+(V_\RR)$ be the spin and vector representations.
The complex structure $I_X$ of $X$ acts on the spin representation $H^*(X,\RR)$ 
via $m(\tilde{I})$, for an element $\tilde{I}\in\Spin(V_\RR)_B$, and $\rho(\tilde{I})\in SO_+(V_\RR)$ is the complex structure of $X\times\hat{X}$. 
The adjoint orbit in $\Spin(V_\RR)_B$ of $\tilde{I}$ is shown to be the period domain for a complete family of polarized abelian varieties of Weil type deformation equivalent to $(X\times\hat{X},\eta)$.

In Section \ref{sec-Examples} we construct a complex multiplication $\eta:K\rightarrow \End_\QQ(X\times\hat{X})$ and a polarization $\Xi_t$, such that $(X\times\hat{X},\eta,\Xi_t)$ is a polarized abelian variety of split Weil type. The construction depends on an embedding $\hat{\eta}:F\rightarrow\End_\QQ(X)$, a
polarization $\Theta\in \wedge^2_FH^1(X,\QQ)$, an element $q\in F$ with $K=F(\sqrt{-q})$, and an element $t\in K_-$.

In Section \ref{sec-orlov-equivalence} we consider two coherent sheaves $F_1$, $F_2$ on $X$ with $ch(F_i)\in B$. We refer to such sheaves as {\em secant sheaves}. We associate to $F_1$ and $F_2$ an object $\G$ in $D^b(X\times\hat{X})$,
with a $\Spin(V)_B$-invariant characteristic class $\kappa(\G):=ch(\G)\exp(-c_1(\G)/\rank(\G))$, via a derived equivalence $\Phi:D^b(X\!\times \!X)\rightarrow D^b(X\!\times\!\hat{X})$ introduced by Orlov.
The $\Spin(V)_B$-invariance of $\kappa(\G)$ implies that the class $\kappa_{d/2}(\G)\in H^{\frac{d}{2},\frac{d}{2}}(X\times\hat{X},\QQ)$ is the sum $\gamma+\delta$, where $\gamma$ is a Weil class and $\delta$ is a polynomial in classes in $\A^2$.
In Section \ref{sec-BB} we give a criterion in terms of $ch(F_i)$, $i=1,2$, for the Weil class $\gamma$ associated to $\G$ not to vanish.

In Section \ref{sec-semiregular-sheaf-over-sixfold} we specialize to the case where $K=\QQ(\sqrt{-q})$ is an imaginary quadratic number field, $q$ a positive integer, and $X$ is the Jacobian of a genus $3$ curve $C$. 
We show that the ideal sheaf of $q+1$ translates of the Abel-Jacobi curve on $X$ is a secant sheaf. This leads to a pair of secant sheaves $F_1$ and $F_2$, for which the dual of $\G$ is a coherent sheaf $\E$ on $X\times\hat{X}$, which equivariantly  satisfies the hypotheses of the Semi-regularity theorem. 
We then use the theorem to prove the algebraicity of the Weil classes on all polarized abelian sixfolds of split Weil type.

In Section \ref{sec-higher-dimension} we express the hope that a speculative stronger version of the Semi-regularity theorem  would lead to the proof of algebraicity of Weil classes on some abelian varieties of dimension $\geq 8$.

%
\section{The Semi-regularity Theorem}
\label{sec-semiregularity}
Let $K$ be a CM-field and $F$ its totally real subfield. 
Let $X$ be an abelian variety admitting an algebra embedding $\hat{\eta}:F\rightarrow \End_\QQ(X)$. 
Set $\hat{X}:=\Pic^0(X)$. The rough strategy is to (1) construct an embedding $\eta:K\rightarrow \End_\QQ(X\times\hat{X})$ making $X\times\hat{X}$ an abelian variety of Weil type. (2) Construct a coherent sheaf $E$ on $X\times\hat{X}$ a characteristic class of which yields an algebraic Weil class in $HW(X\times\hat{X},\eta)$. (3) Deform the sheaf $E$ along with $(X\times\hat{X},\eta)$ to all deformation equivalent abelian varieties of Weil type using the Semi-regularity theorem.

We state next the Semi-regularity theorem. Let $Y$ be an $N$-dimensional compact K\"{a}hler manifold and let $E$ be a coherent sheaf over $Y$.
Denote by $at_E\in\Ext^1(E,E\otimes\Omega^1_Y)$ the Atiyah class of $E$. Let $\sigma_q$, $0\leq q\leq N-2$, be the composition
\[
\Ext^2(E,E)\RightArrowOf{at_E^q/q!} \Ext^{2+q}(E,E\otimes \Omega^q_Y)\RightArrowOf{tr} H^{q,q+2}(Y).
\]
The semi-regularity map is 
\[
\sigma_E:=(\sigma_0, \dots, \sigma_{N-2}):\Ext^2(E,E)\rightarrow \oplus_{q=0}^{N-2} H^{q,q+2}(Y).
\]
The sheaf $E$ is {\em semi-regular}, if $\sigma_E$ is injective. Note that if $Y$ is a surface, then $\sigma_E=\sigma_0$. If, furthermore, $Y$ is a $K3$ or abelian surface and $E$ is a simple sheaf, i.e., $\End(E)\cong\CC$, then $\sigma_0$ is an isomorphism, by Serre's duality. Hence, every simple sheaf over a $K3$ or abelian surface is semi-regular. 

Buchweitz and Flenner show that the following diagram is commutative
\begin{footnotesize}
\begin{equation}
\label{eq-Diagram-semi-regularity}
\xymatrix{
H^1(Y,TY) \ar[rr]^{\Contract at_E} \ar[rd]_{\Contract ch(E)}
&& \Ext^2(E,E) \ar[ld]^{\sigma_E}
\\
& \oplus_{q=0}^{N-2} H^{q,q+2}(Y)
}
\end{equation}
\end{footnotesize}

\noindent
\cite[Cor. 4.3]{buchweitz-flenner}.
A first order infinitesimal deformation $\xi\in H^1(TY)$ of $Y$ belongs to the kernel of $\Contract at_E$, if and only if it can be lifted to a first order deformation of the pair $(Y,E)$ \cite{toda}. The kernel of $\Contract ch(E)$ consists of $\xi\in H^1(TY)$, such that $ch(E)$ remains of Hodge type in the direction of $\xi$, by Griffiths' transversality. If $E$ is semi-regular, the two conditions are equivalent, by the commutativity of the diagram. Buchweitz and Flenner extended this observation, generalizing Bloch's Semi-regularity theorem \cite{bloch} as follows. Let $\pi:\Y\rightarrow B$ be a deformation of a compact K\"{a}hler manifold $Y_0$ over a smooth germ $(B,0)$ and set $Y_b:=\pi^{-1}(b)$ for $b\in B$.  Let $\gamma=\sum_{p=0}^N\gamma_p,$ where $\gamma_p$ is a Hodge class in $H^{p,p}(Y_0,\QQ)$.
We say that $\gamma$ {\em remains of Hodge type over $B$}, if for all $p$ the class $\gamma_p$ extends to a horizontal section of the local system $R^{2p}\pi_*\QQ$, which belongs to the direct summand $R^p\pi_*\Omega^p$ under the Hodge decomposition. 

\begin{thm}\cite[Th. 1.5]{buchweitz-flenner}
\label{thm-semiregularity}
Let $E$ be a semi-regular coherent sheaf over $Y_0$, such that $ch(E)$ remains of Hodge type over $B$.
Then $E$ extends to a coherent sheaf over $\pi^{-1}(U)$ for some open analytic neighborhood $U$ of $0$ in $B$.
\end{thm}

Assume instead that the base $B$ is a smooth and connected analytic space and the fibers of $\pi$ are projective. 
The theorem then implies that each class $ch_p(E)$ remains algebraic over the whole of $B$, since the locus where it is algebraic is a countable union of Zariski closed analytic subsets \cite[Sec. 4.2]{voisin}, and it contains the non-empty open subset $U$ over which $E$ deforms.

Assume that the rank $r$ of $E$ does not vanish.
The condition that $ch(E)$ remains of Hodge type over $B$ is equivalent to the conjunction of two independent conditions:
\begin{enumerate}
\item
\label{cond-c-1}
The class $c_1(E)$ remains of Hodge type over $B$.
\item
\label{cond-kappa}
The class $\kappa(E):=ch(E)\exp(-c_1(E)/r)$  remains of Hodge type over $B$.
\end{enumerate}
Assume that (\ref{cond-kappa}) holds but (\ref{cond-c-1}) fails. 
If $c_1(E)/r$ is integral, then there exists a line bundle $L$ over $Y_0$ with $c_1(L)=c_1(E)/r$ and $ch(E\otimes L^{-1})=\kappa(E)$
remains of Hodge type over $B$. The tensor product of a semi-regular sheaf with a line-bundle is semi-regular and so Theorem \ref{thm-semiregularity} applies to 
$E\otimes L^{-1}$ to conclude that $\kappa(E)$ remains algebraic in every fiber of $\pi$.

The above conclusion remains valid even 
if $c_1(E)/r$ is not integral. 
In that case we need to replace $L$ above by a line bundle over a $\mu_r$-gerb. Equivalently and more elementary, we replace $L$ by a rank $1$ locally free coherent sheaf
twisted by a \v{C}ech $2$-cocycle $\theta$ with coefficient in the local system $\mu_r$ of $r$-th roots of unity. See \cite{calduraru-thesis} for basic facts about twisted sheaves.
Represent the line bundle $\det(E)$ by a \v{C}ech $1$-cocycle 
$\varphi_{ij}\in\StructureSheaf{Y_0}^*(U_{ij})$ with respect to an open covering $\U:=\{U_i\}_{i\in I}$ in the analytic topology, such that $U_{ij}:=U_i\cap U_j$
is simply connected, for all $i,j\in I$. Let $\psi_{ij}$ be an $r$-th root of $\varphi_{ij}$ and set 
$\theta_{ijk}:=\psi_{ij}\psi_{jk}\psi_{ki}^{-1}$. Then $\theta_{ijk}^r=1$, since $\{\varphi_{ij}\}$ is a cocycle, and $\theta:=\{\theta_{ijk}\}$ is a  \v{C}ech $2$-cocycle in
$Z^2(\U,\mu_r)$. The $2$-cocycle $\theta$ is the coboundary associated to the $1$-cochain $\psi_{ij}\in C^1(\U,\StructureSheaf{}^*),$ so the class of $\theta$ is in the kernel of $H^2(\U,\mu_r)\rightarrow H^2(X,\StructureSheaf{}^*)$. 
The $1$-cochain $\{\psi_{ij}\}$ glues $\{\StructureSheaf{U_i}\}_{i\in I}$ to a rank $1$ coherent sheaf $L$ twisted by the cocycle $\theta$.
The tensor product $E\otimes L^{-1}$ is a coherent sheaf twisted by the cocycle $\theta^{-1}$. The line-bundle $\det(E\otimes L^{-1})$ is twisted by the $2$-cocycle $\theta^{-r}$, which is trivial, and is represented by the gluing cocycle $\varphi_{ij}\psi_{ij}^{-r}=1$, and so the determinant line bundle of $E\otimes L^{-1}$ is trivial.
The definitions of the Atiyah class, Chern character, and the semi-regularity map, all extend to $\mu_r$-twisted sheaves (\cite{lieblich},\cite[Def. 7.3.5, 7.3.6]{markman-sixfolds}), the class
$\kappa(E)$ is the Chern character $ch(E\otimes L^{-1})$ \cite[Lem. 7.3.7]{markman-sixfolds}, and $E\otimes L^{-1}$ is semi-regular, since $E$ is by assumption.
The semi-regularity theorem was generalized by Pridham to the setting of derived stacks, which applies to our $\mu_r$-gerb \cite[Rem. 2.26]{pridham} (see also \cite[Sec 7.4]{markman-sixfolds} for a proof in the case of families of abelian varieties).
We thus conclude that if $E$ is semi-regular and $\kappa(E)$ remains of Hodge type over $B$, then $\kappa(E)$ remains algebraic in every fiber of $\pi$. 
%
\section{Abelian varieties and spin groups}
\label{sec-spin-groups}

We recall in this section the fundamental role spin groups play in the geometry of abelian varieties.
%
Let $X$ be an abelian $n$-fold and set $\hat{X}:=\Pic^0(X)$. Then $H^1(\hat{X},\ZZ)$ is naturally isomorphic to $H^1(X,\ZZ)^*$ and so the abelian group 
$V:=H^1(X,\ZZ)\oplus H^1(\hat{X},\ZZ)$ is endowed with the symmetric bilinear pairing 
$((w_1,\theta_1),(w_2,\theta_2))_V:= \theta_1(w_2)+\theta_2(w_1)$,
which is unimodular, even, of signature $(2n,2n)$.
Note that the two direct summands $H^1(X,\ZZ)$ and $H^1(\hat{X},\ZZ)$ are each a maximal isotropic subgroup.
Let $SO(V)$ be the special orthogonal group of $V$ and let $SO(V_\bullet)$ be that of $V_\bullet:=V\otimes_\ZZ\bullet$, where $\bullet$ is any field.
Let $\Spin(V)$ be the subgroup of $\Spin(V_\QQ)$ leaving the lattice $V\subset V_\QQ$ invariant. Let $\rho:\Spin(V_\bullet)\rightarrow SO(V_\bullet)$ be the natural homomorphism and
denote its image by $SO_+(V_\bullet)$, for $\bullet=\ZZ$ or a field. 

The spin representation of $\Spin(V_\bullet)$  is constructed as the exterior algebra of a maximal isotropic subspace of $V_\bullet$ \cite{chevalley}.
We choose the subspace to be $H^1(X,\bullet)$, so that $S_\bullet:=H^*(X,\bullet)$ is the spin representation and $S^+_\bullet:=H^{ev}(X,\bullet)$ and $S^-_\bullet:=H^{odd}(X,\bullet)$ are the half-spin representations. 

Set $S\!\!:=\!\!H^*(X,\!\ZZ)$. It comes with the natural bilinear pairing
$
(\alpha,\!\beta)_S\!\!:=\!\!\int_X\!\tau(\alpha)\cup\beta,
$
where $\tau:H^i(X,\ZZ)\rightarrow H^i(X,\ZZ)$ is multiplication by $(-1)^{i(i-1)/2}$. The pairing is symmetric, if $n$ is even, and anti-symmetric, if $n$ is odd.
It extends naturally to a pairing on $S_\QQ$. Given an object $F$ in the bounded derived category of coherent sheaves $D^b(X)$ we have $\tau(ch(F))=ch(F^\vee)$,
where $F^\vee$ is the derived dual object $\RHom(F,\StructureSheaf{X})$.

The spin representation is defined over $\ZZ$ as follows. Set $Q(v):=(v,v)_V/2$.
The Clifford algebra $C(V)$ is the quotient of the tensor algebra $\oplus_{k=0}^\infty V^{\otimes k}$ by the two-sided ideal generated by 
$\{ v^2-Q(v) \ : \ v\in V\}$, where the integer $Q(v)$ 
is regarded in $V^{\otimes 0}:=\ZZ$. 
The algebra $C(V)$ is $\ZZ/2\ZZ$-graded, $C(V)=C(V)^{ev}\oplus C(V)^{odd}$. 
Let 
$
m:V\rightarrow \End(S)
$
be the homomorphism mapping $(w,\theta)\in V$ to $w\wedge(\bullet)+\theta\Contract(\bullet)$, where $\theta\Contract:S\rightarrow S$ is the contraction with $\theta\in H^1(\hat{X},\ZZ)\cong H^1(X,\ZZ)^*$. Then $m_{v_1}m_{v_2}+m_{v_2}m_{v_1}=(v_1,v_2)id_S$, and so $m$ extends to a homomorphism
\begin{equation}
\label{eq-m-C-V}
m:C(V)\rightarrow \End(S),
\end{equation}
which is in fact an isomorphism of algebras \cite[Prop. 3.2.1(e)]{golyshev-luntz-orlov}. Note that $V$ projects injectively into $C(V)$ and we denote its image by $V$ as well. The Clifford group is the subgroup $G(V)$ of invertible elements of $C(V)$ conjugating $V$ to itself
\[
G(V):=\{g\in C(V)^\times \ : \ gVg^{-1}=V\}
\]
and we denote by 
\begin{equation}
\label{eq-rho-G-V}
\rho:G(V)\rightarrow O(V)
\end{equation} 
the natural homomorphism. Given an element $v\in V$ with $(v,v)=\pm 2$, then $v$ is an element of $G(V)$
and $\rho_{v}$ is minus the reflection in the co-rank $1$ sublattice $v^\perp\subset V$ orthogonal to $v$. The {\em main anti-involution} $\tau:C(V)\rightarrow C(V)$
sends $v_1\cdots v_r$ to $v_r\cdots v_1$. The main involution $\alpha:C(V)\rightarrow C(V)$
acts on $C(V)^{ev}$ as the identity and multiplies $C(V)^{odd}$ by $-1$. The conjugation $x\mapsto x^*$ is the composition $\tau\circ\alpha.$
Then
\[
\Spin(V)=\{g\in C(V)^{ev} \ : \ gg^*=1 \ \mbox{and} \ gVg^{-1}=V
\}.
\]
The spin representation
\begin{equation}
\label{eq-spin-representation-m}
m:\Spin(V)\rightarrow \Aut(S,(\bullet,\bullet)_S)
\end{equation}
is the restriction of (\ref{eq-m-C-V}) and is hence faithful. The vector representation
$
\rho:\Spin(V)\rightarrow SO_+(V)
$
is the restriction of (\ref{eq-rho-G-V}) and its kernel is $\{\pm 1\}$.

Let $\Aut(D^b(X))$ be the group of isomorphism classes of exact auto-equivalences of $D^b(X)$. Every equivalence of derived categories $\Phi:D^b(X)\rightarrow D^b(Y)$ between two smooth projective varieties $X$ and $Y$ is represented as the {\em Fourier-Mukai transform}
\[
R\pi_{Y,*}(L\pi_X^*(\bullet)\stackrel{L}{\otimes} F):D^b(X)\rightarrow D^b(Y),
\]
where $\pi_X$ and $\pi_Y$ are the projections from $X\times Y$, by a theorem of Orlov. The object $F$ in $D^b(X\times Y)$ is unique up to isomorphism and is called the {\em Fourier-Mukai kernel} of $\Phi$. We denote by $\Phi^H:H^*(X,\QQ)\rightarrow H^*(Y,\QQ)$ the correspondence homomorphism $[ch(F)]^*$. Denote by $\Spin_{Hdg}(V)\subset \Spin(V)$ the subgroup preserving the Hodge structure of $V=H^1(X\times\hat{X},\ZZ)$.
Orlov proved that when $X$ and $Y$ are abelian varieties, then $\Phi^H:H^*(X,\ZZ)\rightarrow H^*(Y,\ZZ)$ is an integral isomorphism \cite{orlov-abelian-varieties}.
Furthermore, when $X=Y$, then the homomorphism $\Phi\mapsto \Phi^H$ factors through a surjective homomorphism 
$\Aut(D^b(X))\rightarrow \Spin_{Hdg}(V)$ and the spin representation (\ref{eq-spin-representation-m}). We get the short exact sequence
\begin{equation}
\label{eq-Orlov-short-exact-sequence}
0\rightarrow X\times\hat{X}\times 2\ZZ \rightarrow \Aut(D^b(X))\rightarrow \Spin_{Hdg}(V)\rightarrow 0,
\end{equation}
where the subgroup $X$ corresponds to the auto-equivalences induced by translation automorphisms and the subgroup
$\hat{X}$ corresponds to tensorization by line bundles in $\Pic^0(X)$ (see \cite{orlov-abelian-varieties}). The factor $2\ZZ$ corresponds to even shifts $[2k]$, $k\in\ZZ$. 
The subgroup $X\times\hat{X}$ is the identity component of $\Aut(D^b(X))$ and the latter acts on it by conjugation. The resulting homomorphism
$\Aut(D^b(X))\rightarrow GL(H^1(X\times\hat{X},\ZZ))$ factors through $\Spin_{Hdg}(V)$ via the restriction of the vector representation $\rho$ \cite{orlov-abelian-varieties,golyshev-luntz-orlov}. 

%
\section{A strategy for proving the algebraicity of Weil classes}
\label{sec-strategy}
Let $K$ be a CM-field and let $F$ be its totally real subfield. 
Assume that the abelian $n$-fold $X$ admits an algebra embedding $\hat{\eta}:F\rightarrow\End_\QQ(X)$.
For every rational representation $\psi:F\rightarrow \End(\QQ^{2n})$ (so that $n=de/4$), there exist  non-empty $\frac{n}{2}\left(\frac{d}{2}+1\right)$-dimensional moduli spaces of polarized abelian $n$-folds $X$ with an algebra embedding $\hat{\eta}:F\rightarrow\End_\QQ(X)$ equivalent to $\psi$ \cite[Sec. 9.2]{BL}.
See \cite{ellenberg,shimada} and references therein for examples of Jacobians with real multiplication by $F$.
We get an algebra embedding $F\rightarrow \End_\QQ(X\times\hat{X})$, which we denote by $\hat{\eta}$ as well.
Let $\hat{\Sigma}$ be the set of all field embeddings $F\rightarrow \RR$.
The vector space structure of $V_\QQ=H^1(X\times\hat{X},\QQ)$ over $F$ yields the decomposition 
$V_\RR=\oplus_{\hat{\sigma}\in\hat{\Sigma}}V_{\hat{\sigma},\RR}$, where the pairing $(\bullet,\bullet)_V$ restricts to the natural pairing on $V_{\hat{\sigma},\RR}=H^1_{\hat{\sigma}}(X,\RR)\oplus H^1_{\hat{\sigma}}(\hat{X},\RR)$. 

We would like to extend $\hat{\eta}$ to an embedding $\eta:K\rightarrow\End_\QQ(X\times\hat{X})$ endowing $X\times\hat{X}$ with the structure of an abelian variety of Weil type, 
establish the algebraicity of the Weil classes on $X\times\hat{X}$, 
and use the Semi-regularity theorem to deform them to all abelian varieties of Weil type in the same connected component in moduli. 
Let $\A^2\subset H^{1,1}(X\times\hat{X},\QQ)$ be the $e/2$-dimensional subspace of classes that remain of Hodge type under all deformations of $(X\times\hat{X},\eta)$ as an abelian variety of Weil type. An explicit description of $\A^2$ is given in Proposition \ref{prop-invariant-algebra-A} below.
Following is our strategy.

\begin{enumerate}
\item
\label{strategy-cond1}
Construct a complex multiplication $\eta:K\rightarrow \End_\QQ(X\times\hat{X})$, extending $\hat{\eta}$, and an $\eta$-compatible polarization $h$. Let  
$V_\CC=\oplus_{\sigma\in\Sigma}V_\sigma$ be the decomposition associated to $\eta$.
We require that for each $\sigma:K\rightarrow\CC$, restricting to $\hat{\sigma}:F\rightarrow\RR$, the summand $V_\sigma$ is a maximal isotropic subspace of $V_{\hat{\sigma},\CC}$.
\item
\label{strategy-cond2}
Construct a coherent sheaf $E$ over $X\times\hat{X}$ of non-zero rank $r$ satisfying:
\begin{enumerate}
\item
\label{strategy-cond-semi-regularity}
$E$ is semi-regular.
\item
\label{strategy-cond-remains-Hodge}
The class $\kappa(E):=ch(E)\exp(-c_1(E)/r)$ remains of Hodge type under all deformations of $(X\times\hat{X},\eta,h)$ as a polarized abelian variety of Weil type.
\item
\label{strategy-cond-linear-independence}
The class $\kappa_{d/2}(E)$ in $H^{\frac{d}{2},\frac{d}{2}}(X\times\hat{X},\QQ)$ does not belong to the image of $\Sym^{d/2}(\A^2)$.
\end{enumerate}
\end{enumerate}

Conditions 
(\ref{strategy-cond-semi-regularity}) and 
(\ref{strategy-cond-remains-Hodge}) and 
the Semi-regularity theorem imply that $\kappa(E)$ remains algebraic on every polarized abelian variety of Weil type $(A,\eta',h')$
in the connected component of moduli containing $(X\times\hat{X},\eta,h)$. It follows that $\kappa_{d/2}(E)$ belongs to 
the subspace $Im[\Sym^{d/2}(\A^2)]\oplus HW(X\times\hat{X},\eta)$, since these are the Hodge classes that remain of Hodge type on the generic abelian variety of Weil type, by Equation (\ref{eq-Hodge-classes-of-degree-d-on-the-generic-A}). It follows that $\gamma:=\kappa_{d/2}(E)-\delta$ is a non-zero class in $HW(X\times\hat{X},\eta)$, for some $\delta\in Im[\Sym^{d/2}(\A^2)]$, by
Condition 
(\ref{strategy-cond-linear-independence}). The class $\gamma$ remains algebraic on every $(A,\eta',h')$ deformation equivalent to  $(X\times\hat{X},\eta,h)$, since $\kappa_{d/2}(E)$ and $\delta$ do. Now $K$ acts on $H^*(A,\QQ)$ via algebraic correspondences and $HW(A,\eta')$ is $1$-dimensional over $K$. Hence every class in $HW(A,\eta')$ is algebraic.

Condition 
(\ref{strategy-cond-remains-Hodge}) imposes a compatibility between the class $\kappa_{d/2}(E)$ and the complex multiplication $\eta$. The latter is equivalent to the data of the character subspaces $V_{\sigma,\CC}$, $\sigma\in\Sigma$, of $V_\CC$. The key to achieving this compatibility is Chevalley's theory of pure spinors, which is our next topic. We require the subspaces $V_\sigma$ to be maximal isotropic in Condition (\ref{strategy-cond1}) in order for the representation theory of spin groups to guide our construction.

%
\section{Pure spinors}
\label{sec-pure-spinor}
The Grassmannian $IGr(2n,V_\bullet)$ of maximal isotropic subspaces of $V_\bullet$, $\bullet$ a field, has two connected components, $IGr^+(2n,V_\bullet)$ and $IGr^-(2n,V_\bullet)$. 
Given a maximal isotropic subspace $W\subset V_\bullet$, set
\begin{equation}
\label{eq-ell-W}
\ell_W:=\{
\lambda\in S_\bullet \ : \ m_v(\lambda)=0, \  \forall v\in W
\}.
\end{equation}
Then $\ell_W$ is a one-dimensional subspace of $S_\bullet$, which is contained either in $S^+_\bullet$, if $W\in IGr^+(2n,V_\bullet)$, or in $S^-_\bullet$, if $W\in IGr^-(2n,V_\bullet)$ \cite[III.1.4]{chevalley}. The morphism
\[
\ell: IGr^\pm(2n,V_\bullet)\rightarrow \PP(S^\pm_\bullet)
\]
is a $\Spin(V_\bullet)$-equivariant embedding and its image is called the {\em even/odd spinorial variety}. A non-zero element $\lambda\in \ell_W$ is called an {\em even/odd  pure spinor} and we will refer to $\ell_W$ as a pure spinor as well.
It follows immediately from the definition (\ref{eq-ell-W}) that 
\[
\ell_{H^1(\hat{X},\bullet)}=H^0(X,\bullet) \ \ \ \mbox{and} \ \ \ 
\ell_{H^1(X,\bullet)}= H^{2n}(X,\bullet).
\]
In particular, $1\in H^0(X,\bullet)\subset S^+_\bullet$ is an even pure spinor. Given an element $g\in \Spin(V_\bullet)$, the class $m_g(1)$ is thus an even pure spinor.

\begin{example}
\label{example-pure-spinor}
Let $K$ be a CM-field with a totally real subfield $F$. Choose $q\in F$, such that  $K=F(\sqrt{-q}),$
where $\sqrt{-q}$ is a choice of a square root in $K$, so that $\hat{\sigma}(q)>0$, for all $\hat{\sigma}\in\hat{\Sigma}$. Assume that the abelian $n$-fold $X$ admits an embedding $\hat{\eta}:F\rightarrow \End_\QQ(X)$.
Let $\Theta\in \wedge^2_FH^1(X,\QQ)$ be a non-degenerate class. 
Set
\begin{eqnarray}
\label{eq-pure-spinor}
\exp(\sqrt{-q}\Theta)&:=&1+\sqrt{-q}\Theta + \cdots  +\frac{(\sqrt{-q})^k}{k!}\Theta^k \cdots =\alpha+\sqrt{-q}\beta,
\\
\nonumber
\alpha &:=& 1-\frac{q}{2}\Theta^2 +\cdots + \frac{(-q)^j}{(2j)!}\Theta^{2j}+\cdots
\\
\nonumber
\beta&:=&\Theta-\frac{q}{3!}\Theta^3+\cdots+\frac{(-q)^j}{(2j+1)!}\Theta^{2j+1}+\cdots
\end{eqnarray}
where here $\Theta^k$ denotes the element of $\wedge^{2k}_FH^1(X,\QQ)$. 
We have the natural isomorphism
\[
\Hom_F(H^1(X,\QQ),F)\IsomRightArrow \Hom_\QQ(H^1(X,\QQ),\QQ),
\]
given by $f\mapsto tr_{F/\QQ}\circ f$. Hence, the embedding $\hat{\eta}$ extends to an algebra embedding $\hat{\eta}:F\rightarrow\End_F(V_\QQ)$ and we denote by $V_{\hat{\eta}}$ the vector space $V_\QQ$ regarded as a vector space over $F$.
Denote by 
\begin{equation}
\label{eq-F-valued-pairing-on-V-QQ}
(\bullet,\bullet)_{V_{\hat{\eta}}}
\end{equation}
the natural $F$-valued pairing on $V_{\hat{\eta}}$, so that $(x,y)_V=tr((x,y)_{V_{\hat{\eta}}})$.
We get the groups $\Spin(V_{\hat{\eta}})$ and $\Spin(V_{\hat{\eta}}\otimes_F K)$. Product with $\exp(\sqrt{-q}\Theta)$ in the exterior algebra $(\wedge^*_FV_{\hat{\eta}})\otimes K\cong \wedge^*_K (V_{\hat{\eta}}\otimes_FK)$ is an element of $m(\Spin(V_{\hat{\eta}}\otimes_F K))$.
Evaluating it at the pure spinor $1$ we see that the element $\exp(\sqrt{-q}\Theta)$ is an even pure spinor corresponding to a maximal isotropic 
subspace $W\subset V_{\hat{\eta}}\otimes_FK$. Let $g\in \Spin(V_{\hat{\eta}}\otimes_FK)$ be the element satisfying $m_g=\exp(\sqrt{-q}\Theta)\cup(\bullet)$.
One computes that
\begin{equation}
\label{eq-W-Theta-q}
W=\rho_g(H^1(\hat{X},\QQ)\otimes_FK)=\{(-\sqrt{-q}(\theta\Contract \Theta),\theta) \ : \ 
\theta\in \Hom_F(H^1(X,\QQ),F)_F\otimes K\}.
\end{equation}
Note that $W\cap \iota(W)=(0)$, where $\iota$ acts on the second tensor factor of $V_{\hat{\eta}}\otimes_FK$ as the involution in $\Gal(K/F)$, since $\Theta$ is assumed non-degenerate.
\EndProof
\end{example}

Assume given a maximal isotropic subspace $W\subset V_{\hat{\eta}}\otimes_FK$ such that $W\cap\iota(W)=(0)$.
If $F=\QQ$ and $K=\QQ(\sqrt{-q})$, $q$ a positive integer, then a choice of $\sqrt{-q}\in\CC$ determines an embedding $\sigma:K\rightarrow\CC$ and an embedding of $W$ as a maximal isotropic subspace $V_\sigma$ of 
$V_\CC$. The complex conjugate embedding $\bar{\sigma}$ yields the maximal isotropic subspace $V_{\bar{\sigma}}=\bar{V}_\sigma$.
For a general CM-field $K$ we get a maximal isotropic subspace $W_T$ of $V_\CC$ associated to each choice of a CM-type $T:\hat{\Sigma}\rightarrow\Sigma$.
We denote by $T$ the two equivalent data:
\begin{itemize}
\item
a subset of $\Sigma$ consisting of a choice of one embedding for each pair of two complex conjugate embeddings, and
\item
a right inverse $T:\hat{\Sigma}\rightarrow \Sigma$ of the restriction map $\Sigma\rightarrow \hat{\Sigma}$.
\end{itemize} 
We have the isomorphisms
\begin{equation}
\label{eq-decomposition-of-V-CC}
V_{\hat{\eta}}\otimes_\QQ\CC\cong V_{\hat{\eta}}\otimes_F(F\otimes_\QQ\RR)\otimes_\RR\CC=\oplus_{\hat{\sigma}\in\hat{\Sigma}}(V_{\hat{\eta}}\otimes_{F,\hat{\sigma}}\RR)\otimes_\RR\CC.
\end{equation}
Denote by $W_{T,\hat{\sigma}}$ the subspace spanned over $\CC$ by the image of $W$ via 
the homomorphism 
\begin{equation}
\label{eq-homomorphism-defining-W-T}
id\otimes\hat{\sigma}\otimes T(\hat{\sigma}):V_{\hat{\eta}}\otimes_F F \otimes_FK\rightarrow V_{\hat{\eta}}\otimes_{F,\hat{\sigma}}\RR\otimes_\RR\CC.
\end{equation}
$W_{T,\hat{\sigma}}$ is a subspace of the direct summand corresponding to $\hat{\sigma}$ on the right hand side of (\ref{eq-decomposition-of-V-CC}). Set
\[
W_T:=\oplus_{\hat{\sigma}\in\hat{\Sigma}}W_{T,\hat{\sigma}}.
\]
One checks that $W_T$ is a maximal isotropic subspace of $V_\CC$, for every CM-type $T$. Indeed, $W_{T,\hat{\sigma}}$ is a maximal isotropic subspace of 
$V_{\hat{\eta}}\otimes_{F,\hat{\sigma}}\CC$, for each $\hat{\sigma}\in\hat{\Sigma}$.

The direct sum decomposition $H^1(X,\CC)=\oplus_{\hat{\sigma}\in\hat{\Sigma}} H^1_{\hat{\sigma}}(X,\CC)$ yields the isomorphism
\begin{equation}
\label{eq-tensor-product-factorization-of-S}
S_\CC:=\wedge^* H^1(X,\CC) \cong \otimes_{\hat{\sigma}\in\hat{\Sigma}} \wedge^*H^1_{\hat{\sigma}}(X,\CC),
\end{equation}
where the tensor product is in the category of $\ZZ/2\ZZ$ graded algebras. Denote  the even pure spinor 
$\ell_{W_{T,\hat{\sigma}}}\subset \wedge^{ev}H^1_{\hat{\sigma}}(X,\CC)$ by $\ell_{T,\hat{\sigma}}$. Then $\otimes_{\hat{\sigma}\in\hat{\Sigma}}\ell_{T,\hat{\sigma}}\subset S^+_\CC$ 
is the pure spinor $\ell_{W_T}$. Let $\T_K$ be the set of all CM-types for $K$.
The linear subspace $B_\CC\subset S^+_\CC$ spanned by the lines $\{\ell_{W_T}\}_{T\in\T_K}$ is defined over $\QQ$ and corresponds to a subspace
\begin{equation}
\label{eq-B}
B\subset S^+_\QQ
\end{equation}
of dimension $2^{e/2}$ \cite[Lem. 7.1.3, Cor. 7.2.2]{markman-CM}. 

Denote by $\bar{T}$ the CM-type complex conjugate to $T$ given by  $\bar{T}(\hat{\sigma})=T(\hat{\sigma})\circ\iota.$
Then $W_{\bar{T},\hat{\sigma}}=\overline{W}_{T,\hat{\sigma}}$. Let $P_{\hat{\sigma}}\subset \wedge^{ev}H^1_{\hat{\sigma}}(X,\CC)$ be the $2$-dimensional subspace spanned by $\ell_{T,\hat{\sigma}}$ and $\ell_{\bar{T},\hat{\sigma}}$. We have the equality $B_\CC=\otimes_{\hat{\sigma}\in\hat{\Sigma}} P_{\hat{\sigma}}$ with respect to the factorization (\ref{eq-tensor-product-factorization-of-S}).
%
\section{Complex Multiplication}
\label{sec-complex-multiplication}
Assume given a maximal isotropic subspace $W\subset V_{\hat{\eta}}\otimes_FK$ such that $W\cap\iota(W)=(0)$.
Then $V_\QQ=\{w+\iota(w) \ : \ w\in W\}$.
Define
$
\eta: K\rightarrow \End_F(V_{\hat{\eta}}\otimes_FK)
$
by letting $\eta_t$ act on $W$ by multiplication by $t$ and on $\iota(W)$ by multiplication by $\iota(t)$, for all $t\in K$. 
Then for $w\in W$, $\eta_t(w+\iota(w))=tw+\iota(tw)$ and so $\eta_t(V_\QQ)=V_\QQ$. We get the algebra embedding
\begin{equation}
\label{eq-eta}
\eta:K\rightarrow \End(V_\QQ).
\end{equation}

The subspace $W_{T,\hat{\sigma}}$ of $V_{\hat{\sigma},\RR}$ depends only on the value $T(\hat{\sigma})$, as the homomorphism 
(\ref{eq-homomorphism-defining-W-T}) depends only on this value. Hence, given $t\in K$, $\eta_t$ acts on $W_{T,\hat{\sigma}}$ via multiplication by 
$T(\hat{\sigma})(t)$. So the subspace 
$
V_\sigma
$
 of $V_\CC$ on which $\eta_t$ acts via $\sigma(t)$ is $W_{T,\hat{\sigma}}$, if and only if $T(\hat{\sigma})=\sigma$ and
 \begin{equation}
 \label{eq-K-character-decomposition-of-W-T}
 W_{T,\CC}=\oplus_{\hat{\sigma}\in\hat{\Sigma}}V_{T(\hat{\sigma})}.
 \end{equation}

Let $\Spin(V_\QQ)_B$ be the subgroup of $\Spin(V_\QQ)$ fixing every point in the secant space $B\subset S^+_\QQ$ associated to $W$ in (\ref{eq-B}).
Define $\Spin(V_\RR)_B$ analogously. Note that $\rho$ maps $\Spin(V_\QQ)_B$ injectively into $SO_+(V_\QQ)$, since $-1$ does not belong to $\Spin(V_\QQ)_B$. 

Let $\Spin(V_\bullet)_\eta$ be the subgroup of $\Spin(V_\bullet)$ consisting of elements commuting with $\eta(K)$, for $\bullet=\QQ$ or $\RR$. 

\begin{lem}
\label{lemma-Spin-V-B-commutes-with-eta-K}
$\Spin(V_\bullet)_B$ is a subgroup of 
$\Spin(V_\bullet)_\eta$, for $\bullet=\QQ$ or $\RR$. 
\end{lem}

\begin{proof}
Let $g$ be an element of $\Spin(V_\bullet)_B$. Each $\ell_{W_T}$ is
$m_g$-invariant, for all $T\in \T_K$. Thus,  $W_T$ is $\rho_g$-invariant, for all $T\in \T_K$. 
Let $\sigma\in\Sigma$ restrict to $F$ as $\hat{\sigma}$. 
The subspace $V_\sigma$ is the intersection of $W_T$ and $W_{T'}$, if $T(\hat{\sigma})=T'(\hat{\sigma})=\sigma$ and $T(\hat{\sigma}')\neq T'(\hat{\sigma}')$,
for all $\hat{\sigma}'\in\hat{\Sigma}\setminus\{\hat{\sigma}\}$. Hence, $V_\sigma$ is $\rho_g$-invariant, for all $\sigma\in\Sigma$. Consequently, $\rho_g$ commutes with $\eta(K)$.  
\end{proof}

\begin{lem}
\label{lemma-B-Hodge-implies-eta-is}
\cite[Lem. 7.4.1]{markman-CM}
Assume that $B$ is spanned by Hodge classes. 
Then $\eta(K)$ is contained in $\End_{Hdg}(V_\QQ)\cong\End_\QQ(X\times\hat{X})$.
Furthermore, $\dim(V^{1,0})_\sigma=\dim(V^{0,1})_\sigma$, for all $\sigma\in\Sigma$. Consequently, $(X\times\hat{X},\eta)$ is of Weil type.
\end{lem}

\begin{proof}[Sketch of proof that $\eta(K)\subset \End_{Hdg}(V_\QQ)$] 
There exists an element $\tilde{I}$ of $\Spin(V_\RR)$ satisfying $m_{\tilde{I}}=I_X$ and $\rho_{\tilde{I}}=I_{X\times\hat{X}}.$
The assumption that $B$ is spanned by Hodge classes implies that $\tilde{I}$ belongs to $\Spin(V_\RR)_B$. 
Lemma \ref{lemma-Spin-V-B-commutes-with-eta-K} implies that 
$\rho_{\tilde{I}}$ commutes with $\eta(K)$. Hence, 
$\eta(K)$ is contained in $\End_{Hdg}(V_\QQ)$. 
\end{proof}

%

%
\section{The $\Spin(V_\QQ)_B$-invariant subalgebra of $H^*(X\times\hat{X},\QQ)$}
\label{sec-A}
 Let $K_-$ be the $-1$ eigenspace of $\iota:K\rightarrow K$.

\begin{lem}
\label{lemma-H-t}
\begin{enumerate}
\item
\cite[Lem 5.1.2]{markman-CM}
The endomorphism $\eta_{\iota(t)}$ is the adjoint of $\eta_t$ with respect to $(\bullet,\bullet)_V$, for all $t\in K$.
\item
\label{lemma-item-Xi-t-in-Spin-V-eta-invariant}
\cite[Cor. 5.2.1]{markman-CM}
Given $t\in K_-$, the $F$-valued 
pairing $\Xi_t(x,y):=(\eta_t(x),y)_{V_{\hat{\eta}}}$, where $(x,y)_{V_{\hat{\eta}}}$ is given in (\ref{eq-F-valued-pairing-on-V-QQ}), 
is anti-symmetric and $\Spin(V_\QQ)_\eta$-invariant. 
The resulting homomorphism\footnote{
Here we again regard $\wedge^2_FV^*_\QQ$ as a subspace of $\wedge^2_\QQ V^*_\QQ$.
The subspace $\oplus_{\hat{\sigma}\in\hat{\Sigma}}\wedge^2 V^*_{\hat{\sigma},\RR}$ of $\wedge^2V^*_\RR$ is defined over $\QQ$ and corresponds to the image of the injective homomorphism  
$\wedge^2_FV^*_\QQ\rightarrow \wedge^2_\QQ V^*_\QQ$ given by $\theta(\bullet,\bullet)\mapsto tr_{F/\QQ}\circ \theta(\bullet,\bullet).$
} 
$\Xi:K_-\rightarrow \wedge^2_FV^*_\QQ\cap (\wedge^2V_\QQ^*)^{\Spin(V_\QQ)_\eta}$ is injective.
\item
\label{lemma-item-H-is-Spin-V-B-invariant}
\cite[Lem. 7.3.1]{markman-CM}
Given $t\in K_-$, let $H_t:V_\QQ\times V_\QQ\rightarrow K$ be given by
\[
H_t(x,y):=(-t^2)(x,y)_{V_{\hat{\eta}}}+t\Xi_t(x,y).
\]
Then $H_t$ is a $\Spin(V_\QQ)_B$-invariant hermitian form. 
\end{enumerate}
\end{lem}

Given $g\in SO_+(V_{\hat{\sigma},\RR})$ and $\sigma\in\Sigma$ restricting to $\hat{\sigma}\in\hat{\Sigma}$, $g$ leaves the restriction of $\sigma\circ H_t$ invariant $\Leftrightarrow$ $g$ leaves the restriction of $\sigma\circ \Xi_t$ invariant $\Leftrightarrow$ $(\eta_t(x),y)_V=(\eta_t(g(x)),g(y))_V$, for all $x,y\in V_{\hat{\sigma},\RR}$, $\Leftrightarrow$
$g$ commutes with the restriction of $\eta_t$ to $V_{\hat{\sigma},\RR}$. Hence, if $g$ leaves the restriction of $\sigma\circ H_t$ invariant, then $g$ leaves $V_\sigma$ and $V_{\bar{\sigma}}$ invariant. Let 
$SU(V_{\hat{\sigma},\RR})$ be the subgroup of $SO_+(V_{\hat{\sigma},\RR})$ of elements leaving the restriction of $\sigma\circ H_t$ invariant and restricting
to each of $V_{\sigma,\CC}$ and $V_{\bar{\sigma},\CC}$ with determinant $1$. 

\begin{lem}
\label{lemma-SU}
\cite[Lem 9.1.1]{markman-CM}
The natural homomorphism
$
\Spin(V_\RR)_B\rightarrow \prod_{\hat{\sigma}\in\hat{\Sigma}}SO_+(V_{\hat{\sigma},\RR})
$
maps $\Spin(V_\RR)_B$ isomorphically onto $\prod_{\hat{\sigma}\in\hat{\Sigma}}SU(V_{\hat{\sigma},\RR})$.
\end{lem}

Lemma \ref{lemma-Spin-V-B-commutes-with-eta-K} shows that $V_\sigma$ is $\Spin(V_\QQ)_B$-invariant  and 
Lemma \ref{lemma-SU} shows that $\rho_g$, $g\in \Spin(V_\QQ)_B,$ acts trivially on $\wedge^dV_\sigma$, for all $\sigma\in\Sigma$.
Hence, $HW(X\times\hat{X},\eta)$ consists of  $\Spin(V_\QQ)_B$-invariant classes. 
Lemmas \ref{lemma-H-t}(\ref{lemma-item-Xi-t-in-Spin-V-eta-invariant}) and \ref{lemma-Spin-V-B-commutes-with-eta-K} show that $\Xi(K_-)$ is an $e/2$-dimensional subspace of $\Spin(V_\QQ)_B$-invariant elements. 
Let $\A:=(\wedge^*V_\QQ)^{\Spin(V_\QQ)_B}=H^*(X\times\hat{X},\QQ)^{\Spin(V_\QQ)_B}$ be the $\Spin(V_\QQ)_B$-invariant subalgebra. 

\begin{prop}
\label{prop-invariant-algebra-A}
\begin{enumerate}
\item
\label{prop-item-generators-for-A}
\cite[Prop. 8.0.1]{markman-CM}
The subalgebra $\A$ is generated by $\A^2$ and 
$HW(X\times\hat{X},\eta)$. The graded summand $\A^2$ is the image of $\Xi(K_-)$ in $\wedge^2V_\QQ$ via the isomorphism $\wedge^2V_\QQ^*\cong \wedge^2V_\QQ$ induced by the pairing $(\bullet,\bullet)_V$.
\item
\label{prop-item-A-is-Hodge}
If $B$ is spanned by Hodge classes, then so is $\A$.
\end{enumerate}
\end{prop}

\begin{proof}[Proof of (\ref{prop-item-A-is-Hodge})] Assume that $B$ is spanned by Hodge classes.
The subspace $HW(X\times\hat{X},\eta)$ consists of Hodge classes and  $I_{X\times\hat{X}}$ belongs to $\rho(\Spin(V_\RR)_B)$ and commutes with $\eta(K)$, by Lemma \ref{lemma-B-Hodge-implies-eta-is}. Now $I:=I_{X\times\hat{X}}$ is an isometry with respect to $(\bullet,\bullet)_V$. Hence,
\[
\Xi_t(I(x),I(y)):=(\eta_t(I(x)),I(y))_V=(I\eta_t(x),I(y))_V=(\eta_t(x),y)_V=\Xi_t(x,y)
\] 
and so $\Xi_t$ is of type $(1,1)$. Hence, $\A$ consists of Hodge classes, by part (\ref{prop-item-generators-for-A}).
\end{proof}

\begin{lem}
\label{lemma-adjoint-orbit-is-a-period-space}
Assume that $B$ is spanned by Hodge classes and $\Xi_t\in H^{1,1}(X\times\hat{X},\QQ)$ is an ample class, for some $t\in K_-$. 
Then the adjoint orbit of $I_{X\times\hat{X}}$ in $\Spin(V_\RR)_B$ consists of complex structures $I$ on the differentiable manifold $A$ underlying $X\times\hat{X}$
with respect to which 
\begin{enumerate}
\item  \cite[Lem. 9.2.2, 9.2.3]{markman-CM}
$((A,I),\eta,\Xi_t)$ is a polarized abelian variety of Weil type.
\item
\label{lemma-item-A-consists-of-Hodge-classes}
The subalgebra $\A:=H^*(A,\QQ)^{\Spin(V_\QQ)_B}$ of \ $\Spin(V_\QQ)_B$-invariant classes in $H^*(A,\QQ)$ consists of Hodge classes.
\end{enumerate}
\end{lem}

\begin{proof}[Proof of (\ref{lemma-item-A-consists-of-Hodge-classes})]
Note that the proof of Proposition \ref{prop-invariant-algebra-A}(\ref{prop-item-A-is-Hodge}) applies to any element $I$ of $\Spin(V_\RR)_B$ satisfying $I^2=-id$.
\end{proof}

An example of complex multiplication $\eta$ on $X\times\hat{X}$ satisfying the hypothesis of the above Lemma is given in the next section.
The adjoint orbit ${\mathcal O}$ of $I_{X\times\hat{X}}$ in $\Spin(V_\RR)_B$ parametrizes a complete family of polarized abelian varieties of Weil type,
by \cite[Rem. 9.2.4]{markman-CM}.
Let $T$ be a CM-type and consider the map ${\mathcal O}\rightarrow \prod_{\hat{\sigma}\in \hat{\Sigma}}Gr(d/2,V_{T(\hat{\sigma})})$
sending a complex structure $I$ to $\{V_{T(\hat{\sigma})}^{1,0}\}_{\hat{\sigma}\in\hat{\Sigma}}$.
The above map is an embedding of ${\mathcal O}$ as an open subset, in the classical topology, of the product of  grassmannians \cite[Sec. 9.2]{markman-CM}. See also \cite{deligne-milne} and \cite[Lemma 11.5.25]{charles-schnell}.

The group $\Spin(V_\QQ)_B$ is the special Mumford-Tate group\footnote{
The special Mumford-Tate group of $H^1((A,I),\QQ)$ is the smallest algebraic subgroup of $GL(H^1(A,\RR))$, which is defined over $\QQ$ and which contains the circle group $\{a+bI  :  a,b\in\RR, \ a^2+b^2\!=\!1\}$.} 
of the generic abelian variety of Weil type 
deformation equivalent to $(X\times\hat{X},\eta)$.
%

%
\section{Examples of polarized $X\times\hat{X}$ of split Weil type}
\label{sec-Examples}
Assume that the class $\Theta\in\wedge^2_FH^1(X,\QQ)$ in Example \ref{example-pure-spinor} is such that $tr_{F/\QQ}\circ \Theta(\bullet,\bullet)$ is an ample class in $H^{1,1}(X,\QQ)$. If $X$ is simple and $\hat{\eta}:F\rightarrow \End_\QQ(X)$ is an isomorphism, then $H^{1,1}(X,\QQ)\subset \wedge^2_FH^1(X,\QQ)$
(see \cite[Lem. 11.1.1]{markman-CM}) and so any polarization $\Theta$ would arise this way. Let $W\subset V_{\hat{\eta}}\otimes_FK$ be the maximal isotropic subspace associated to $(\Theta, \sqrt{-q})$ in (\ref{eq-W-Theta-q}). 
Let $\eta:K\rightarrow \End_{Hdg}H^1(X\times\hat{X},\QQ)$ be the embedding given in (\ref{eq-eta})
and let $B\subset H^{ev}(X,\QQ)$ be the secant linear subspace given in (\ref{eq-B}), both associated to $W$. 
Let $t\in K_-$ be a non-zero element. 

\begin{lem} 
\label{lemma-X-times-hat-X-of-Weil-type}
\cite[Lemma 11.1.2]{markman-CM}
\begin{enumerate}
\item
\label{lemma-item-B-is-spanned-by-Hodge-classes}
$\eta(K)\subset \End_{Hdg}H^1(X\times\hat{X},\QQ)$, 
\item
\label{lemma-item-of-split-type}
$H_t$ is of split-type, and 
\item
\label{lemma-item-ample}
$t$ can be chosen so that $\Xi_t$ is an ample class and $(X\times\hat{X},\eta,\Xi_t)$ is a polarized abelian variety of split Weil type.
\end{enumerate}
\end{lem}

\begin{proof}[Sketch of proof]
(\ref{lemma-item-B-is-spanned-by-Hodge-classes}) It suffices to prove that $B$ is spanned by Hodge classes, by Lemma \ref{lemma-B-Hodge-implies-eta-is}. We know that $B$ is a rational subspace, so it suffices to prove that it is contained in $\oplus_{p=0}^{2n}H^{p,p}(X)$. The proof follows easily from the assumption that $\Theta$ is of type $(1,1)$.

(\ref{lemma-item-of-split-type})
Let $\{y_1,\dots, y_d\}$ be an $F$-basis of $H^1(\hat{X},\QQ)$, such that $\span\{y_1, \dots, y_{d/2}\}$ is a $\Theta$-isotropic subspace. Let $g\in \Spin(V_{\hat{\eta}}\otimes_FK)$
be the element in Example \ref{example-pure-spinor} satisfying the equality $m_g(\bullet)=\exp(\sqrt{-q}\Theta)\cup(\bullet)$.
Then $\{\rho_g(y_1), \dots, \rho_g(y_d)\}$ is a $K$-basis of $W$, by definition of $W$, and so
$
\{\rho_g(y_j)+(id_{V_{\hat{\eta}}}\otimes\iota)(\rho_g(y_j)) \ : \ 1\leq j\leq d\}
$
is a $K$-basis of $V_\QQ$ with respect to $\eta$. We have seen in Equation (\ref{eq-W-Theta-q}) that 
$\rho_g(y_j)=(-\sqrt{-q}(y_j\Contract \Theta),y_j)$. Hence,
\[
\rho_g(y_j)+(id_{V_{\hat{\eta}}}\otimes\iota)(\rho_g(y_j))=2(0,y_j)
\]
and $\{(0,y_1),\dots, (0,y_d)\}$ is a $K$-basis of $V_\QQ$. The definition of $\eta$ yields the equality 
$\eta_t(0,y)=(-\hat{\eta}_{t\sqrt{-q}}(y\Contract\Theta),0)$.
We claim that the subspace $Z:=\span_K\{(0,y_j)\}_{j=1}^{d/2}$
is $H_t$ isotropic. Indeed, it is $(\bullet,\bullet)_{V_{\hat{\eta}}}$-isotropic, since $H^1(\hat{X},\QQ)$ is $(\bullet,\bullet)_{V_{\hat{\eta}}}$-isotropic, 
and it is $\Xi_t$ isotropic, since 
\begin{eqnarray*}
\Xi_t((0,y_j),(0,y_k))\!&\!\!=\!\!&\!
(\eta_t(0,y_j),(0,y_k))_V=
(\hat{\eta}_{t\sqrt{-q}}(-y_j\Contract\Theta,0),(0,y_k))_V=-t\sqrt{-q}\Theta(y_j,y_k)
\\
&\!\!=\!\!&0, \ \mbox{for} \ 1\leq j,k\leq d/2.
\end{eqnarray*}

(\ref{lemma-item-ample}) 
We regard $tr_{F/\QQ}\circ\Xi_t\in \wedge^2_\QQ V_\QQ^*$ as an element of $\wedge^2_\RR V_\RR^*$ by extending scalars.
The symmetric bilinear form $(tr_{F/\QQ}\circ\Xi_t)(\bullet,I_{X\times\hat{X}}(\bullet))$ restricts to a definite form on $V_{\hat{\sigma},\RR}$, for each $\hat{\sigma}\in\hat{\Sigma}$,
by the argument of \cite[Lem. 3.1.1]{markman-sixfolds}. The choice of $t$ affects the sign of this form. 
The image of the map $F\rightarrow \RR^{\hat{\Sigma}}$, given by $f\mapsto (\hat{\sigma}(f))_{\hat{\sigma}\in\hat{\Sigma}}$, contains a full lattice.  Hence, for every function $s:\hat{\Sigma}\rightarrow\{+,-\}$ there exists $f\in F$, such that the sign of $\hat{\sigma}(f)$ is $s(\hat{\sigma})$, for all $\hat{\sigma}\in\hat{\Sigma}$.
It follows that $t\in K_-$ can be chosen, so that the bilinear pairing is positive definite on $V_{\hat{\sigma},\RR}$, 
 for each $\hat{\sigma}\in\hat{\Sigma}$. The equality $\dim(V^{1,0}_\sigma)=\dim(V^{0,1}_\sigma)=\frac{d}{2}$, for all $\sigma\in\Sigma$, follows from Lemma \ref{lemma-B-Hodge-implies-eta-is}.
\end{proof}

{\bf Strategy revisited:}
Note that Lemma \ref{lemma-X-times-hat-X-of-Weil-type} achieves item (\ref{strategy-cond1}) in the Strategy section \ref{sec-strategy}. 
The compatibility condition (\ref{strategy-cond-remains-Hodge}) in Strategy section \ref{sec-strategy}, between $\kappa(E)$ and $\eta$, is equivalent to the condition that $\kappa(E)$ would be $\Spin(V_\QQ)_B$-invariant, by Proposition \ref{prop-invariant-algebra-A}(\ref{prop-item-generators-for-A}) and
Lemma \ref{lemma-adjoint-orbit-is-a-period-space}(\ref{lemma-item-A-consists-of-Hodge-classes}). 
If $F_1$ and $F_2$ are coherent sheaves on $X$ with Chern characters $ch(F_i)\in B$, then
$F_1\boxtimes F_2:=\pi_1^*F_1\otimes\pi_2^*F_2$ is a 
coherent sheaf on $X\times X$ with $ch(F_1\boxtimes F_2)$ invariant with respect to the diagonal action of $\Spin(V_\QQ)_B$. We would like to transform such a sheaf to an object $E$ on $X\times\hat{X}$ with a $\Spin(V_\QQ)_B$-invariant $\kappa(E)$ via an equivalence $\Phi:D^b(X\times X)\rightarrow D^b(X\times\hat{X})$. This requires 
the induced isomorphism of cohomologies $\Phi^H:H^*(X\times X)\rightarrow H^*(X\times\hat{X})$ to be $\Spin(V)$-equivariant, up to cup product with $exp(\alpha)$
for some $\alpha\in H^{1,1}(X\times\hat{X},\QQ)$. In the next section we observe that  Orlov introduced such an equivalence in \cite{orlov-abelian-varieties}.
%
\section{Orlov's derived equivalence $\Phi:D^b(X\times X)\rightarrow D^b(X\times\hat{X})$}
\label{sec-orlov-equivalence}
Let $\mu:X\times X\rightarrow X\times X$ be the automorphism given by $\mu(x,y)=(x+y,y)$. 
Let $\P$ be the Poincar\'{e} line bundle over $\hat{X}\times X$, normalized so that it restricts trivially to $\hat{X}\times\{0\}$. 
Let $id\times\Phi_\P:D^b(X\times\hat{X})\rightarrow D^b(X\times X)$  be the integral transform with Fourier-Mukai kernel 
$\pi_{13}^*\StructureSheaf{\Delta_X}\otimes\pi_{24}^*\P$, where $\pi_{ij}$ is the projection from $X\times \hat{X}\times X\times X$ to the product of the $i$-th and $j$-th factors. Orlov's derived equivalence
$
\Phi:D^b(X\times X)\rightarrow D^b(X\times\hat{X})
$
is the inverse of $\mu_*\circ (id\times\Phi_\P):D^b(X\times\hat{X})\rightarrow D^b(X\times X)$. 
The Fourier-Mukai kernel $\U\in D^b(X\times \hat{X}\times X\times X)$ of $\mu_*\circ (id\times\Phi_\P)$ restricts to $\{(x,L)\}\times X\times X$
as the Fourier-Mukai kernel $\U_{x,L}\in D^b(X\times X)$ of the auto-equivalence $L\otimes t_{x,*}:D^b(X)\rightarrow D^b(X)$ of push-forward via the translation
automorphism $t_x(y)=x+y$ of $X$ and tensorization by the line bundle $L$. So $\Phi^{-1}$ encodes the fact that $X\times\hat{X}$ is a subgroup of the group of auto-equivalences of $D^b(X)$.

Let $\Phi^H:H^*(X\times X,\ZZ)\rightarrow H^*(X\times\hat{X},\ZZ)$
be the correspondence induces by the Chern character of the Fourier-Mukai kernel of $\Phi$. Let $\tau:H^i(X,\ZZ)\rightarrow H^i(X,\ZZ)$ be multiplication by $(-1)^{i(i-1)/2}$. 
The K\"{u}nneth theorem identifies $H^*(X\times X,\ZZ)$ as the tensor square $S\otimes_\ZZ S$ of the spin representation and $H^*(X\times\hat{X},\ZZ)$
is the exterior algebra $\wedge^*V$ of the vector representation. Hence, both are endowed with the structure of an integral representation of $\Spin(V)$. The two are not isomorphic, but they become isomorphic\footnote{The representations are reducible, so the set of isomorphisms is not a $\QQ^{\times}$-torsor.} 
after tensoring with $\QQ$ \cite[Sec. III.3.3]{chevalley}.

\begin{prop}
\label{prop-Spiv-V-equivariance-of-Orlov-equivalence}
\cite[Prop. 6.1.2]{markman-sixfolds}
The following composition is $\Spin(V)$-equivariant
\[
\tilde{\phi}:=\exp(-c_1(\P)/2)\cup \Phi^H\circ (id\otimes \tau): H^*(X\times X,\QQ)\rightarrow H^*(X\times\hat{X},\QQ).
\]
\end{prop}

The integral isomorphism 
\[
\Phi^H\circ (id\otimes \tau): S\otimes_\ZZ S\cong H^*(X\times X,\ZZ)\rightarrow H^*(X\times\hat{X},\ZZ)\cong\wedge^*V
\]
conjugates the diagonal action of $m_g$, $g\in \Spin(V)$, on $S\otimes_\ZZ S$ to an automorphism of $\wedge^*V$ leaving invariant the decreasing  filtration
$F_k(\wedge^*V):=\oplus_{j\geq k}\wedge^jV$ and it was previously known that the induced action on the graded summand $\wedge^kV$ is $\wedge^k\rho_g$
\cite[Prop. 4.3.7 and Cor. 4.3.8]{golyshev-luntz-orlov}. If $g$ belongs to $\Spin_{Hdg}(V)$, then the action of $m_g$ on $H^*(X\times X)$ need not preserve the grading and is the diagonal cohomological action of an auto-equivalence of $X$, by the right exactness of Orlov's sequence (\ref{eq-Orlov-short-exact-sequence}), while the action of $\rho_g$ is the cohomological action of an automorphism of $X\times\hat{X}$. The paper \cite{golyshev-luntz-orlov} refers to this correspondence of symmetries as an instance of homological Mirror-Symmetry.

Keep the assumptions and notation of Section \ref{sec-complex-multiplication}. Let $F_1$ and $F_2$ be coherent sheaves on $X$ with $ch(F_i)$ in the secant linear subspace $B\subset H^{ev}(X,\QQ)$.
We say that $F_1$ and $F_2$ are {\em secant sheaves}. Given a class $\gamma$ in $H^*(X\times\hat{X},\QQ)$ with graded summand $\gamma_i$ in $H^{2i}(X\times\hat{X},\QQ)$ and with non-zero $\gamma_0\in \QQ$,
set $\kappa(\gamma):=\gamma\exp(-\gamma_1/\gamma_0)$. Given an object $E\in D^b(\bullet)$ of non-zero rank, we have $\kappa(E):=\kappa(ch(E))$.

\begin{cor}
\label{cor-kappa-E-is-Spin-V-B-invariant}
Assume that the rank of the object $E:=\Phi(F_1\boxtimes F_2^{\vee})\in D^b(X\times\hat{X})$ is non-zero. 
The class $\kappa(E)$ is $\Spin(V_\QQ)_B$-invariant.
\end{cor}

\begin{proof}
The class $\gamma\!:=\!\tilde{\phi}(ch(F_1\boxtimes F_2))\!=\!\exp(-c_1(\P)/2)\!\cup\! ch(E)$ is $\Spin(V)_B$-invariant, by the $\Spin(V)$-equivariance of $\tilde{\phi}$ in Proposition \ref{prop-Spiv-V-equivariance-of-Orlov-equivalence} and the assumed $\Spin(V_\QQ)_B$-invariance of $ch(F_1\boxtimes F_2)$.
It follows that $\gamma_1$ is $\Spin(V_\QQ)_B$-invariant and hence so is $\kappa(\gamma)$. Finally, $\kappa(\gamma)=\kappa(E)$.
\end{proof}

%
\section{The homomorphism $B\otimes B\rightarrow HW(X\times\hat{X},\eta)$}
\label{sec-BB}
Keep the notation of Corollary \ref{cor-kappa-E-is-Spin-V-B-invariant}.
We consider next Condition (\ref{strategy-cond-linear-independence}) in the Strategy section \ref{sec-strategy}.
It requires the class $\kappa_{d/2}(E)\in \A^d\oplus HW(X\times\hat{X},\eta)$ not to belong to the direct summand $\A^d$.
We need to translate it to a condition on $ch(F_1)$ and $ch(F_2)$.

We define a natural grading $B\otimes_\QQ B=\oplus_{k=0}^{e/2} BB_k$ on
$B\otimes_\QQ B$ as follows. 
Given two CM-types $T, T'$, set $T\cap T':=T(\hat{\Sigma})\cap T'(\hat{\Sigma})\subset\Sigma$
and let $|T\cap T'|$ be its cardinality. The subspace 
\[
BB_{k,\CC}:= \bigoplus_{(T_1,T_2)\in \T_K\times\T_K, \ |T_1\cap T_2|=k } \ell_{W_{T_1}}\otimes \ell_{W_{T_2}}
\]
is defined over $\QQ$ and corresponds to a subspace $BB_k\subset B\otimes_\QQ B$ of dimension $\Choose{e/2}{k}2^{e/2}$.
So $\dim(BB_1)=e2^{(\frac{e}{2}-1)}$.  The equality 
$W_T\cap W_{T'}=\oplus_{\hat{\sigma}\in T\cap T'}V_{T(\hat{\sigma})}$ holds, by Equation (\ref{eq-K-character-decomposition-of-W-T}).

\begin{lem}
\cite[Lemma 10.1.3]{markman-CM} 
If $|T\cap T'|=k$, then the line $\Phi^H(\ell_{W_T}\otimes \tau(\ell_{W_{T'}}))$ belongs to $F_{dk}(V_\CC):=\oplus_{i\geq dk}(\wedge^iV_\CC)$ and it projects onto the line $\wedge^{dk}[W_T\cap W_{T'}]$ in $\wedge^{dk}V_\CC$. In particular, the subspace $HW(X\times\hat{X},\eta)$ is equal to the image of the composition
\begin{equation}
\label{eq-BB-1-to-HW}
BB_1\LongRightArrowOf{\Phi^H\circ(id\otimes\tau)} F_d(\wedge^*V_\QQ)\rightarrow \wedge^dV_\QQ.
\end{equation}
\end{lem}

\begin{proof}[Sketch of proof]
The key is the analogous statement \cite[III.3.3]{chevalley} relating the top exterior power  $\wedge^{\dim(W+ W')}(W+ W')$ to the image of the tensor product $\ell_W\otimes \ell_{W'}$ of the two pure spinors via 
an isomorphism $\tilde{\varphi}:S\otimes_\ZZ S\rightarrow \wedge^*V$ Chevalley constructs in terms of the Clifford algebra $C(V)$.
The equality
$\Phi^H\circ(id\otimes\tau)=[\Phi_\P^H\otimes(\Phi_\P^H)^{-1}]\circ\tilde{\varphi}$  is proved in \cite[Lem. 6.1.1]{markman-sixfolds}.
 Let $\varsigma:X\times\hat{X}\rightarrow \hat{X}\times X$ be the transposition of the factors.
The isomorphism $[\Phi_\P^H\otimes(\Phi_\P^H)^{-1}]\circ\varsigma_*:\wedge^*V\rightarrow\wedge^*V$ is an analogue of the Hodge star operator and it interchanges $\wedge^iV$ and $\wedge^{4n-i}V$
\cite[Rem. 6.3.3]{markman-sixfolds}.
\end{proof}

We have  $\dim(HW(X\times\hat{X},\eta))=e$. Hence,
$\dim(BB_1)=\dim(HW(X\times\hat{X},\eta))$ and (\ref{eq-BB-1-to-HW}) is an isomorphism, if and only if $e=2$, if and only if $F=\QQ$. In that case $B_\CC=\ell_W\oplus\ell_{\bar{W}}$ 
and Condition (\ref{strategy-cond-linear-independence}) in section \ref{sec-strategy} requires that $ch(F_1)\otimes ch(F_2)$ does not belong to the $2$-dimensional subspace 
$BB_0:=[\ell_W\otimes \ell_{\bar{W}}]\oplus [\ell_{\bar{W}}\otimes\ell_W]$ of the $4$-dimensional $B\otimes_\QQ B$. Classes in $BB_0$ are invariant under the larger group $\Spin(V_\QQ)_\eta$, while classes in $BB_1$ are invariant only under the special Mumford-Tate group $\Spin(V_\QQ)_B$ and so it suffices to check that $ch(F_1)\otimes ch(F_2)$ is not invariant under $\Spin(V_\QQ)_\eta$ (see, for example,  the proof of \cite[Lemma 8.3.1]{markman-sixfolds}).

More generally, write $ch(F_1)\otimes ch(F_2)=\sum_{k=0}^{e/2}\gamma_k$, with $\gamma_k\in BB_k$, and write
$\gamma_1=\sum_{(T,T') \ : \ |T\cap T'|=1}\gamma_{T,T'}$, with $\gamma_{T,T'}\in\ell_{W_T}\otimes\ell_{W_{T'}}$.
Let $\Pi$ be the composition (\ref{eq-BB-1-to-HW}).

\begin{prop}\cite[Prop. 10.2.1(3)]{markman-CM}
Assume that $d>2$, the rank of $\Phi(F_1\boxtimes F_2^\vee)$ is non-zero, and for some $\sigma\in\Sigma$
\begin{equation}
\label{eq-sum-of-gamma-T-T'}
\sum_{(T,T') \ : \ T\cap T'=\{\sigma\}}\Pi( \gamma_{T,T'})\neq 0.
\end{equation} 
 Then $\kappa_{d/2}(E)$
does not belong to the image of $\Sym^{d/2}\A^2$ in $H^{d,d}(X\times\hat{X},\QQ)$.
\end{prop}

\hide{
\begin{example}
Consider the case where $K=F(\sqrt{-q})$, where $q$ is a positive {\em rational} number. 
Let $\tilde{K}$ and $\tilde{F}$ be the Galois closures of $K$ and $F$ over $\QQ$. Then $\Gal(\tilde{K}/\QQ)\cong \Gal(\tilde{F}/\QQ)\times \{id,\iota\}.$
Let $T$ be the CM-type $\{\sigma\in\Sigma \ : \ Im(\sigma(\sqrt{-q})>0\}$. Choose $F_1$ such that $ch(F_1)$ belongs to the plane $P:=\ell_{W_T}\oplus \ell_{W_{\iota(T)}}$, which is defined over $\QQ$. The sum (\ref{eq-sum-of-gamma-T-T'}) consists of at most one non-zero summand, which does not vanish if and only if $ch(F_2)$

\end{example}
}
In the special case when $K$ contains the subfield $\QQ(\sqrt{-q})$, for some positive $q\in\QQ$, then 
the sheaf $F_1$ can be chosen with $ch(F_1)$ in the plane $P$ spanned by $\alpha$ and $\beta$ in Equation (\ref{eq-pure-spinor}).
In that case the sum (\ref{eq-sum-of-gamma-T-T'}) consists of at most one non-zero summand (see \cite[Lem. 11.2.9]{markman-CM}). 
Hence, for such $F_1$ the sum (\ref{eq-sum-of-gamma-T-T'}) does not vanish, for some $\sigma\in\Sigma$,  if and only if the direct summand $\gamma_1$ of $ch(F_1)\otimes ch(F_2)$ in $BB_1$ does not vanish.
An explicit example of coherent sheaves $F_1$ and $F_2$ on the Jacobian of a genus $4$ curve satisfying the hypothesis of the above Proposition  is given in \cite[Ex. 11.2.7 and Lem. 11.2.8]{markman-CM}.

%
\section{A semi-regular sheaf over a sixfold $(X\times\hat{X},\eta)$ of split Weil type}
\label{sec-semiregular-sheaf-over-sixfold}
We finally deal with the semi-regularity condition \ref{strategy-cond-semi-regularity} in the Strategy section \ref{sec-strategy}, but we are currently able to handle it only for $\dim(X)\leq 3$, which forces $K=\QQ(\sqrt{-q})$ to be a quadratic imaginary number field. We may choose $q$ to be a positive integer.
We may also choose $q$ to be even and $\geq 4$, since $\QQ(\sqrt{-q})=\QQ(\sqrt{-4q})$. 

%
\subsection{Two secant sheaves}
\label{sec-two-secant-sheaves}
Let $C$ be a non hyperelliptic curve of genus $3$. Set $X=\Pic^0(C)$. Let $AJ:C\rightarrow \Pic^1(C)$ be the Abel-Jacobi map.
Let $G_1$ and $G_2$ be two cyclic subgroups of $X$ of order $q+1$ satisfying $G_1\cap G_2=(0)$.
Let $C_i\subset X$, $1\leq i\leq q+1$, be $q+1$ disjoint translates of $AJ(C)$, which are transitively permuted by translations by elements of $G_1$. 
Let $\Theta\subset X$ be a translate of the theta divisor in $\Pic^2(C)$. 
Set $F_1:=\Ideal{\cup_{i=1}^{q+1}C_i}(\Theta)$, the ideal sheaf of the union tensored with $\StructureSheaf{X}(\Theta)$. The class Poincar\'{e} dual to $C_i$ is $\Theta^2/2\in H^{2,2}(X,\ZZ)$ and 
one checks that $ch(F_1)=(1-\frac{q}{2}\Theta^2)+(\Theta-\frac{q}{3!}\Theta^3)$,
which belongs to the secant $B=\span_\QQ\{\alpha,\beta\}$, $\alpha=1-\frac{q}{2}\Theta^2$, $\beta=\Theta-\frac{q}{3!}\Theta^3$, associated to the pure spinor
$\exp(\sqrt{-q}\Theta)$ is Example \ref{example-pure-spinor} (see \cite[Lem 8.2.1]{markman-sixfolds} for the computation of $ch(F_1)$).
We get the embedding $\eta:K\rightarrow \End_\QQ(X\times\hat{X})$ and a polarization $\Xi_t$, for some $t\in K_-$, such that $(X\times\hat{X},\eta,\Xi_t)$ is a polarized abelian sixfold of split Weil type, by Lemma \ref{lemma-X-times-hat-X-of-Weil-type}.

Set $C':=-AJ(C)\subset \Pic^{-1}(C)$. Let $C'_i\subset X$, $1\leq i\leq q+1$, be $q+1$ disjoint translates of $C'$, which are transitively permuted by $G_2$. Set $F_2=\Ideal{\cup_{i=1}^{q+1}C'_i}(\Theta)$. 
Then $ch(F_2)=ch(F_1)$ and so $ch(F_2)$ is in $B$ as well. Note that $ch(F_i^\vee)=\alpha-\beta$ is in $B$ as well, where $F_i^\vee$ is the derived dual object.

%
\subsection{The transform of $F_2\boxtimes F_1$ to $X\times\hat{X}$}
If a translate $\tau_t(\Sigma_j)$ of $\Sigma_j$ by a point $-t\in X$ intersects $C_i$, then the two intersect along a length $2$ subscheme and the canonical line bundle of $C_i\cup \tau_{-t}(\Sigma_j)$ is the restriction of 
$L_{i,j,t}(2\Theta)$ for a unique line bundle $L_{i,j,t}\in \hat{X}$ 
(see \cite[Lem. 9.2.7]{markman-sixfolds} for an explicit computation of $L_{i,j,t}$).
Define the morphism $f_{ij}:C_i\times\Sigma_j\rightarrow X\times\hat{X}$  by 
$f_{ij}(x,y)=(y-x,L_{i,j,x-y})$. Denote by $\tilde{\Theta}_{ij}\subset X\times\hat{X}$ the image of $f_{ij}$ and set
$\tilde{\Theta}:=\cup_{1\leq i,j\leq q+1}\tilde{\Theta}_{ij}$. Each $\tilde{\Theta}_{ij}$ projects isomorphically onto a translate of the divisor $\Theta$ in $X$, and $\Theta$ is isomorphic to the symmetric product $C^{(2)}$, by the assumption that $C$ is not hyperelliptic. The genericity assumptions 
\cite[Assum. 9.1.1, 9.2.1]{markman-sixfolds}
assure that the surfaces $\tilde{\Theta}_{ij}$ are pairwise disjoint, and so $\tilde{\Theta}$ is smooth. Let $\Phi$ be Orlov's equivalence. 

\begin{prop}
\label{prop-sheaf-E}
\cite[Lem. 9.1.4, 9.3.1 and Prop. 9.2.2]{markman-sixfolds}
For a suitable choice of the groups $G_1$ and $G_2$ and for a generic $C$, 
the cohomology sheaves $\G_i$ in degree $i$ of the object $\G:=\Phi(F_2\boxtimes F_1)[-3]$ over $X\times\hat{X}$ satisfy:
\begin{enumerate}
\item
$\G_i$ vanishes for $i\not\in\{1,2\}$.
\item
$\G_1$ is a reflexive sheaf of rank $8q$, which is locally free away from $\tilde{\Theta}$. 
\item
$\G_2$ is supported, set theoretically, over $\tilde{\Theta}$.
\item
\label{prop-item-sheaf-E}
The object $\G^\vee[-1]$ is represented by a coherent sheaf $\E$, which is isomorphic to $\SheafHom(\G_1,\StructureSheaf{X\times\hat{X}})$ and $\G_2$ is isomorphic to $\SheafExt^1(\E,\StructureSheaf{X\times\hat{X}})$.
\end{enumerate}
\end{prop}

The fiber of $\G_1$, at a point $(x,L)\in [X\times\hat{X}]\setminus\tilde{\Theta}$, is naturally isomorphic to
\[
H^1(X,\Ideal{\cup_{i=1}^{q+1}C_i}\otimes\Ideal{\cup_{i=1}^{q+1}\tau_{-x}(C'_i)}(\Theta+\tau_{-x}(\Theta))\otimes L^{-1}).
\]

%
\subsection{Equivariance of the object $\G\in D^b(X\times\hat{X})$}
The object 
$\G[3]:=\Phi(F_2\boxtimes F_1)$ is the image of the $G_2\times G_1$-equivariant sheaf $\Ideal{\cup_{i=1}^{q+1}C'_i}\boxtimes \Ideal{\cup_{i=1}^{q+1}C_i}$
via the equivalence $\tilde{\Phi}:=\Phi\circ (\StructureSheaf{X\times\hat{X}}(\Theta\boxtimes\Theta)\otimes(\bullet))$.
Denote by $G_2\times G_1$ also its image in the identity component of $\Aut(D^b(X\times X))$ and let $G:=\tilde{\Phi}\circ (G_2\times G_1)\circ\tilde{\Phi}^{-1}$
be its conjugate in the identity component $(X\times\hat{X})\times\Pic^0(X\times\hat{X})$ of $\Aut(D^b(X\times\hat{X}))$. The equivalence $\tilde{\Phi}$ transforms the natural linearization of $\Ideal{\cup_{i=1}^{q+1}C'_i}\boxtimes \Ideal{\cup_{i=1}^{q+1}C_i}$ to a linearization $\tilde{\lambda}:=\{\tilde{\lambda}_g\}_{g\in G}$ of the object $\G$ with respect to the action of $G$ on $D^b(X\times\hat{X})$, so that $(\G,\tilde{\lambda})$ is an object in the equivariant category $D_G^b(X\times\hat{X})$ in the sense of \cite{BO}.
Given $x\in X$, set $L_x:=\StructureSheaf{X}(\Theta-\tau_x(\Theta))$ and let $\P_x$ be the restriction of the Poincar\'{e} line bundle $\P$ to $\hat{X}\times\{x\}$.
An explicit calculation yields
\[
\tilde{\Phi}\circ(\tau_{x_1},\tau_{x_2})_*\circ\tilde{\Phi}^{-1}=((\pi_1^*L_{x_1}\otimes\pi_2^*\P_{-x_2})\otimes)\circ (\tau_{x_1-x_2},\tau_{L_{x_1+x_2}})_*
\]
(see \cite[Eq. (9.3.1)]{markman-sixfolds}).
We conclude that $G$ projects injectively onto a subgroup $\bar{G}$ of translations in $X\times\hat{X}$, since $G_1\cap G_2=(0)$. 

Let $D:=\det(\G)$ be the determinant line bundle of $\G$.
Let $a$ be a positive integer, such that $8qa\equiv -1$ (mod. $q+1$). Such an integer $a$ exists, by our assumption that $q$ is even.

\begin{lem} \cite[Lem. 9.3.5]{markman-sixfolds}
The object $\G\otimes D^a$ is $\bar{G}$-equivariant. It admits $\bar{G}$-linearization isomorphisms
$\lambda_g:\G\otimes D^a\rightarrow \tau_g^*(\G\otimes D^a)$, $g\in\bar{G}\subset X\times\hat{X},$ satisfying the axioms of a linearization in \cite[Sec. 1.1]{ploog}.
\end{lem}

Let $\pi:X\times\hat{X}\rightarrow A:=(X\times\hat{X})/\bar{G}$ be the quotient morphism onto the quotient abelian variety.  
Let $\bar{\G}$ be the object in $D^b(A)$, which is mapped to $(\G,\lambda)\in D^b_{\bar{G}}(X\times\hat{X})$ via the natural equivalence 
$\pi^*:D^b(A)\rightarrow D^b_{\bar{G}}(X\times\hat{X})$.  Then $\bar{\G}^\vee[-1]$ is represented by a coherent sheaf $\bar{\E}$, by
Proposition \ref{prop-sheaf-E}. We have the isomorphisms
\begin{eqnarray}
\label{eq-Ext-2-bar-E-bar-E}
\Ext^2(\bar{\E},\bar{\E})&\cong& \Hom(\StructureSheaf{A},\bar{\E}^\vee\otimes\bar{\E}[2]) \cong
\Hom(\StructureSheaf{A},\bar{\G}\otimes\bar{\E}[3])
\\
\nonumber
&\cong& 
\Hom(((\StructureSheaf{X\times\hat{X}},1),(\G\otimes D^a,\lambda)\otimes(\G\otimes D^a,\lambda)^\vee[2])^{\bar{G}}
\\
\nonumber
&\cong& \Ext^2(\G\otimes D^a,\lambda),\G\otimes D^a,\lambda))^{\bar{G}}\cong  \Ext^2((\G,\tilde{\lambda}),(\G,\tilde{\lambda}))^G
\\
\nonumber
&\cong& 
\Ext^2(\Ideal{\cup_{i=1}^{q+1}C'_i}\boxtimes\Ideal{\cup_{i=1}^{q+1}C_i},\Ideal{\cup_{i=1}^{q+1}C'_i}\boxtimes\Ideal{\cup_{i=1}^{q+1}C_i})^{G_2\times G_1}
\end{eqnarray}

%
\subsection{The sheaf $\E$ is equivariantly semi-regular}
The second Hochschild cohomology $HH^2(Y)$ of a smooth projective variety $Y$ consists of natural transformations from the identity endo-functor $id$ of $D^b(Y)$ to $id[2]$. Given a coherent sheaf $E$ over $Y$, let $ev_E:HH^2(Y)\rightarrow \Hom(E,E[2])$ be the evaluation 
of a natural transformation on the object $E$.  
The Chern character $ch(E)$ corresponds to a class in the Hochschild homology $HH_*(Y)$, 
via the Hochchild-Kostant-Rosenberg (HKR) isomorphism \cite[Theorem 4.5]{caldararu-II}, and $HH_*(Y)$
is a module over the ring $HH^*(Y)$. 
There is a sufficient criterion for $E$ to be semi-regular \cite[Rem. 8.3.11]{markman-sixfolds}, which we state for $Y$  an abelian variety for simplicity:

\begin{lem}
\label{lemma-sufficient-condition-for-semi-regularity}
\cite[Lem. 8.3.10]{markman-sixfolds}
If the kernel of $ev_E:HH^2(Y)\rightarrow \Hom(E,E[2])$ is equal to the annihilator of $ch(E)$ in $HH^2(Y)$ and $ev_E$ is surjective,
then $E$ is semi-regular.
\end{lem}

\begin{proof}
We have the HKR isomorphism
\[
HH^2(Y)\cong HT^2(Y):=H^2(Y,\StructureSheaf{Y})\oplus H^1(Y,TY)\oplus H^0(Y,\wedge^2TY).
\]
Diagram (\ref{eq-Diagram-semi-regularity}) has an analogous commutative diagram
 replacing 
 $H^1(Y,TY)$ by the first order deformations 
$HH^2(Y)$ 
of the category of coherent sheaves on $Y$ \cite[Prop. 6.2.1 and Cor. 6.3.2]{buchweitz-flenner-HH}. 
The top horizontal homomorphism $\Contract at_E$ in diagram (\ref{eq-Diagram-semi-regularity}) is replaced  
$ev_E:HH^2(Y)\rightarrow \Hom(E,E[2])$. The right arrow  in the new version of diagram (\ref{eq-Diagram-semi-regularity}) is still the semi-regularity map $\sigma_E$ and 
when $Y$ is an abelian variety the left arrow  is the composition of $HKR:HH^2(Y)\rightarrow HT^2(Y)$ with contraction with $ch(E)$. 
\begin{footnotesize}
\[
\xymatrix{
HH^2(Y)\ar[r]^{ev_E}\ar[d]_{HKR} & 
\Ext^2(E,E)\ar[d]^{\sigma_E}
\\
HT^2(E) \ar[r]_-{\Contract ch(E)} & \oplus_{q=0}^{N-2}H^{q,q+2}(Y)
}
\]
\end{footnotesize}
If $\ker(\Contract ch(E)\circ HKR)$ is equal to $\ker(ev_E)$ and $ev_E$ is surjective,
then the equality $\sigma_E\circ ev_E=\Contract ch(E)\circ HKR$ implies that $\sigma_E$ is injective.
\end{proof}

\begin{question}
\label{question-stronger-version-of-semi-regularity}
Is the surjectivity of $ev_E$ needed in Lemma \ref{lemma-sufficient-condition-for-semi-regularity}? Does the 
Semi-regularity Theorem \ref{thm-semiregularity} hold under the weaker assumption that
$\sigma_E$ restricts to the image of $ev_E$ as an injective map (considering also twisted sheaves as we did in Section \ref{sec-semiregularity})?
\end{question}

The sheaf $\E$ of Proposition \ref{prop-sheaf-E}(\ref{prop-item-sheaf-E}) cannot be semi-regular for all $q\geq 4$, since $\dim\Ext^2(\E,\E)$ grows quadratically with $q$, by \cite[Lem. 8.3.8]{markman-sixfolds}, while the co-domain of $\sigma_\E$ is independent of $q$.
Fortunately, $\sigma_\E$ restricts as an injective map to $\Ext^2(\E,\E)^{\bar{G}}$, which contains the image of $ev_\E$. 
Descending to the quotient abelian variety $A:=(X\times\hat{X})/\bar{G}$ enables us to use the semi-regularity theorem and avoid 
Question \ref{question-stronger-version-of-semi-regularity}. 



\begin{prop}
\label{prop-E-is-equivariantly-semi-regular}
The sheaf $\bar{\E}$ over $A$ is semi-regular.
\end{prop}

\begin{proof}[Sketch of proof]
We prove first that the sheaf $\Ideal{C_1}$ over $X$ is semi-regular following
\cite[Lem. 8.3.7]{markman-sixfolds}. 
The surjectivity of $ev_{\Ideal{C_1}}$ follows from the fact that  $HH^1(X)\rightarrow \Ext^1(\Ideal{C_1},\Ideal{C_1})$ is an isomorphism, $HH^1(X)$ generates $HH^*(X)$, $\Ext^1(\Ideal{C_1},\Ideal{C_1})$ generates $\Ext^*(\Ideal{C_1},\Ideal{C_1})$, and so $ev_{\Ideal{C_1}}^*:HH^*(X)\rightarrow\Ext^*(\Ideal{C_1},\Ideal{C_1})$ is a surjective algebra homomorphism. Serre's duality yields $\Ext^2(\Ideal{C_1},\Ideal{C_1})\cong \Ext^1(\Ideal{C_1},\Ideal{C_1})^*$.
It follows that the rank of $ev_{\Ideal{C_1}}:HH^2(X)\rightarrow\Ext^2(\Ideal{C_1},\Ideal{C_1})$ is $6$. One checks that the rank of $\Contract ch(\Ideal{C_1})$ is $6$ as well. The inclusion $\ker(ev_{\Ideal{C_1}})\subset \ker(\Contract ch(\Ideal{C_1}))$ holds, by \cite[Th. B]{Huang}. Hence, $\ker(ev_{\Ideal{C_1}})= \ker(\Contract ch(\Ideal{C_1}))$ and
$\Ideal{C_1}$ is semi-regular, by Lemma \ref{lemma-sufficient-condition-for-semi-regularity}.
The curve $C_1$ maps isomophically onto its image  $\bar{C}$ in $X/G_1$. 
The sheaf $\Ideal{\bar{C}}$ is semi-regular as well. 

Note that the sheaf 
$\Ideal{\cup_{i=1}^{q+1}C_i}$ is naturally $G_1$-equivariant and is the pullback of $\Ideal{\bar{C}}$. 
The group $G_1$ acts trivially on $HH^2(X)$, 
and so $ev_{\Ideal{\cup_{i=1}^{q+1}C_i}}:HH^2(X)\rightarrow \Ext^2(\Ideal{\cup_{i=1}^{q+1}C_i},\Ideal{\cup_{i=1}^{q+1}C_i})^{G_1}$ is surjective and its kernel 
is equal to the kernel of  $\Contract ch(\Ideal{\cup_{i=1}^{q+1}C_i})$ \cite[Cor. 8.3.12]{markman-sixfolds}.
The analogous statement holds for $\Ideal{\cup_{i=1}^{q+1}C'_i}.$ It follows, by the K\"{u}nneth theorem, that 
$ev_{\Ideal{\cup_{i=1}^{q+1}C'_i}\boxtimes \Ideal{\cup_{i=1}^{q+1}C_i}}:HH^2(X\times X)\rightarrow \Ext^2(\Ideal{\cup_{i=1}^{q+1}C'_i}\boxtimes \Ideal{\cup_{i=1}^{q+1}C_i},\Ideal{\cup_{i=1}^{q+1}C'_i}\boxtimes \Ideal{\cup_{i=1}^{q+1}C_i})^{G_2\times G_1}$ is surjective and its kernel is equal to that of
$\Contract ch(\Ideal{\cup_{i=1}^{q+1}C'_i}\boxtimes \Ideal{\cup_{i=1}^{q+1}C_i})$ \cite[Lem. 8.4.1(1) and Lem. 9.3.2]{markman-sixfolds}.
Hence, $ev_{\G}:HH^2(X\times\hat{X})\rightarrow \Ext^2((\G,\tilde{\lambda}),(\G,\tilde{\lambda}))^G$ is surjective and its kernel is equal to that of $\Contract ch(\G)$. 
The analogous statement follows for the sheaf $\E$, by  \cite[Lem. 8.4.1(2)]{markman-sixfolds}, and for the sheaf $\bar{\E}$ by (\ref{eq-Ext-2-bar-E-bar-E}). Hence, $\bar{\E}$ is semi-regular, by Lemma \ref{lemma-sufficient-condition-for-semi-regularity}.
\end{proof}

%
\subsection{Sketch of proof of Theorem \ref{main-thm}}
Step 1 (Weil classes on sixfolds of split Weil type):
We have established 
the semi-regularity of $\bar{\E}$ in Proposition \ref{prop-E-is-equivariantly-semi-regular}. 
The abelian variety $A$ is isogenous to $X\times \hat{X}$ and is thus endowed with the embedding $\eta:K\rightarrow \End_{Hdg}(H^1(X\times\hat{X}))\cong\End_{Hdg}(H^1(A,\QQ))$ and the polarization $\Xi_t\in H^{1,1}(X\times\hat{X},\QQ)\cong H^{1,1}(A,\QQ)$. This completes the implementation of the strategy for the proof of algebraicity of the Weil classes on all deformations of 
$(A,\eta,\Xi_t)$, and hence also for all deformations of $(X\times\hat{X},\eta,\Xi_t)$, as explained in the strategy section \ref{sec-strategy}.
Two connected components of the moduli space of polarized abelian varieties of Weil type of dimension $2n$, the same imaginary quadratic number field, and the same discriminant, parametrize isogenous abelian varieties \cite[Th. 5.2(3)]{van-Geemen}. 
The discriminant takes values in  $\QQ^\times/Nm_{K/\QQ}(K^\times)$ and is the coset of $(-1)^n$ if and only if the component parametrizes polarized abelian varieties of split Weil type \cite[Cor. 4.2]{deligne-milne}.
We conclude that the Hodge Weil classes on all polarized abelian sixfolds of split Weil type are algebraic. 

Step 2 (Weil classes on fourfolds):
The discriminant invariant of polarized abelian varieties with complex multiplication by the same field is multiplicative under cartesian products. Every value in $\QQ^\times/Nm_{K/\QQ}(K^\times)$ is realized  as the discriminant  by some connected component of moduli in every even dimension \cite[Th. 5.2]{van-Geemen}. Hence, for every polarized abelian fourfold $(A_1,\eta_1,h_1)$ of Weil type, of arbitrary discriminant, there exists a polarized abelian surface of Weil type $(A_2,\eta_2,h_2)$, such that the discriminant of their product polarized abelian sixfold of Weil type $(A_1\!\times\! A_2,\eta,\pi_1^*h_1\!+\!\pi_2^*h_2)$ 
is the coset of $-1$. The sixfold is hence of split type and so its Weil classes are algebraic. It follows that the Weil classes of $(A_1,\eta_1,h_1)$ are algebraic, by \cite[Prop. 10]{schoen}.
Hence, the Weil classes are algebraic on every abelian fourfold.
\EndProof

%
\section{What about Weil classes on abelian varieties of dimension $\geq 8$?}
\label{sec-higher-dimension}
We expect that an affirmative answer to Question \ref{question-stronger-version-of-semi-regularity} would lead to a proof of the algebraicity of  Weil classes on some higher dimensional abelian varieties, as well as for CM-fields $K$ with $[K:\QQ]>2$, using the strategy outlined in Section \ref{sec-strategy}. 
It would be interesting to find a systematic construction of examples of secant sheaves on abelian varieties with real multiplication by a totally real field $F$ satisfying  the weaker criterion proposed in 
Question \ref{question-stronger-version-of-semi-regularity}. 
Consider for example a principally polarized abelian fourfold $(X,\Theta)$. Let $d$ be an odd integer $\geq 3$, set $K:=\QQ(\sqrt{-d})$, and set $n:=(d+9)/2$. Let $\{D_i \ : i\in\ZZ/n\ZZ\}$ be $n$ generic translates of the divisor $\Theta$, cyclically indexed, such that for every subset $S\subset \ZZ/n\ZZ$ of cardinality $|S|$ the intersection $\cap_{i\in S}D_i$ is smooth of codimension $|S|$, if $2\leq |S|\leq 4$, or empty if $|S|\geq 5$. Set $Z_i:=D_i\cap D_{i+1}$, $1\leq i\leq n$, and set $Z:=\cup_{i=1}^n Z_i$. Note that the intersection $Z_i\cap Z_j$ is a smooth curve, if $i-j=\pm1$ and it consists of $24$ points if $i\not\in\{j-1,j,j+1\}$.
Let $\nu:\tilde{Z}\rightarrow Z$ be the partial normalization of $Z$ along $(d+9)(2d-1)$ of its $(d+9)(3d-3)$ isolated\footnote{The intersection points of $Z_i\cap Z_j$ are isolated, if and only if $i\not\in\{j-2,j-1,j,j+1,j+2\}$.
} points of self intersection. Then the Chern character of the object 
\[
[\StructureSheaf{X}\RightArrowOf{\nu^*}(\nu_*\StructureSheaf{\tilde{Z}})]\otimes\StructureSheaf{X}(\Theta)
\]
belongs to the secant $\span\{\exp(\sqrt{-d}\Theta),\exp(-\sqrt{-d}\Theta)\}$. 
It is yet to be checked if these secant objects satisfy the weaker criterion in Question \ref{question-stronger-version-of-semi-regularity}.
%
\hspace{0ex}

{\bf Acknowledgements:}
I thank Salvatore Floccari for several helpful comments. I thank Pierre Deligne for insightful comments and corrections. 
This work was partially supported by a grant from the Simons Foundation Travel Support for Mathematicians program (\#962242).


\end{document}